\documentclass{article}

\usepackage{amsmath, amsfonts}
\usepackage[left=1in, right=1in]{geometry}
\usepackage{hyperref}
\usepackage{enumerate}
\usepackage{graphicx}

\newcommand{\given}{\mbox{{\bf $\ \mid \ $}}}

\DeclareMathOperator*{\argmin}{arg\,min}

\newcommand{\slim} {\mathop{\rm lim\,sup}}
\newcommand{\ilim} {\mathop{\rm lim\,inf}}
\newcommand{\field}[1]{\mathbb{#1}}
\newcommand{\PR}{\field{P}}             
\newcommand{\E}{\field{E}}              
\newcommand{\N}{\field{N}}                                
\newcommand{\R}{\field{R}}                                
\newcommand{\A}{\field{A}}                                
\newcommand{\X}{\field{X}}                                
\renewcommand{\S}{{\boldsymbol S}}                                    
\newcommand{\Gr}{\rm Gr}

\newcommand{\K}{\field{K}}
\newcommand{\U}{\field{U}}
\newcommand{\Sp}{\field{S}}
\newcommand{\n}{t}
\renewcommand{\a}{\alpha}
\newcommand{\bS}{\textbf{S}}
\newcommand{\q}{P}
\newcommand{\h}{1}
\newcommand{\B}{\mathcal{B}}
\newcommand{\Y}{\mathbb{Y}}
\newcommand{\T}{\mathbb{T}}
\renewcommand{\P}{\mathbb{P}}
\renewcommand{\c}{\bar c}

\newtheorem{thm}{Theorem}[section]
\newtheorem{cor}[thm]{Corollary}
\newtheorem{lem}[thm]{Lemma}

\newtheorem{defn}[thm]{Definition}
\newtheorem{rem}[thm]{Remark}

\begin{document}

\title{\textbf{Optimality Conditions for Inventory Control}}

\author{\textbf{Eugene A. Feinberg}\medskip \\
  \small{Department of Applied Mathematics and Statistics} \\
  \small{Stony Brook University, Stony Brook, NY 11794-3600, USA}}

\date{}

\maketitle

\begin{abstract}
  This tutorial describes recently developed general optimality conditions for Markov Decision Processes that have significant applications to inventory control.  In particular, these conditions imply the validity of optimality equations and inequalities.  They also imply the convergence of value iteration algorithms.   For total discounted-cost problems only two mild conditions on the continuity of transition probabilities and lower semi-continuity of one-step costs are needed.  For average-cost problems, a single additional assumption on the finiteness of relative values is required.   The general results are applied to periodic-review inventory control problems with discounted and average-cost criteria without any assumptions on demand distributions. The case of partially observable states is also discussed.
\end{abstract}


\section{Introduction}
This tutorial describes  recent progress in the theory of Markov Decision Processes (MDPs) with infinite state and action sets that have significant applications to inventory control.  Two groups of results are covered: (i) optimality conditions for MDPs with total, discounted and average-cost criteria, and (ii) optimality conditions for Partially Observable Markov Decision Processes (POMDPs) with total and discounted cost criteria.

Inventory control studies and applications are important motivating factors for studies of MDPs.  The MDP studies provided important tools for the analysis of inventory control problems.  The parallel development of these fields since the beginning of the second half of the 20th century is broadly recognized.   For example,  the abstract of the historical essay by Girlich and Chikan~\cite{GC} on the history of inventory control studies states:
 ``... we report how inventory problems have motivated the improvement of mathematical disciplines such as Markovian decision theory and optimal control of  stochastic systems to provide a new basis of inventory theory in the second half of our century.''  However, over a long period of time there was a gap between the modeling needs for inventory control, that require mathematical methods for the analysis of infinite-state controlled stochastic systems with unbounded action sets and weakly continuous transition probabilities, and available results for the corresponding models for MDPs.  This gap was recently closed.   Another topic covered in this tutorial is the recent progress in the development of optimality conditions for POMDPs.  The literature on MDPs and inventory control is huge, and we do not attempt a comprehensive survey in this tutorial.  For the most part only directly relevant references are provided.  The reader may find coverage of these topics in the books \cite{BR, besh96, DY, FS, heysol, HLerma1, HL:96, HL:99, Pu, sen99} on MDPs and \cite{Bens0, heysol, Port, SLCB, Zip} on inventory management.

Optimality results for MDPs provide sufficient conditions for the existence of stationary and Markov optimal policies satisfying optimality equations and inequalities, describe continuity properties of the value function, and guarantee the convergence of value iteration and optimal actions when the horizon length tends to infinity or the discount factor tends to 1.  These results provide useful tools to analyze specific inventory control problems and to prove the optimality of particular  policies.  In Section~\ref{s4} this is illustrated with the classic periodic-review single-product stochastic inventory problem with nonnegative  arbitrarily distributed  iid demand.
  Most of the literature on inventory control is limited to discrete or continuous demand distributions.

 Consider the classic periodic-review single-product stochastic inventory problem with backorders.  For a finite horizon and continuous demand,  Scarf~\cite{scarf} established under some conditions the optimality of $(s,S)$ policies. Zabel~\cite{Za} indicated some gaps in \cite{scarf}, corrected them, and mentioned in the last paragraph of \cite{Za} that the proofs there can be adapted to arbitrary demand distributions.  Iglehart~\cite{igle63} and  Veinott and Wagner~\cite{vei65} established the optimality of $(s,S)$ policies for the infinite horizon for continuous and discrete demand respectively.  Zheng~\cite{zhe91} provided an alternative proof for discrete demand.  Beyer and Sethi~\cite{bey99} described and corrected gaps in the proofs in \cite{igle63, vei65}.  As shown in Heyman and Sobel~[Section 7.1]\cite{heysol}, under  appropriate conditions  $(s,S)$ policies are optimal for a finite-horizon problem with arbitrarily distributed demand.  In general, $(s,S)$ policies may not be optimal for finite horizons.  For example, for a problem with convex holding costs the appropriate condition is Assumption~{\bf GB} in Section~\ref{s4}.  This assumption means that, as the amount of backordered inventory increases, the backordering cost per unit time becomes larger than the value of the backordered inventory. However, as shown in Veinott~\cite{vei66a} for discrete demand, $(s,S)$ policies are always optimal for the following three criteria: (i) infinite-horizon average costs per unit time, (ii) infinite-horizon discounted problems with a large discount factor, and (iii) finite-horizon problems with  appropriately selected terminal costs. Chen and Simchi-Levi~\cite{che04b, che04a} described optimal policies for coordinating inventory control and pricing for finite and infinite-horizon problems with general demand under a technical assumption.  If the price is fixed, the
 problem in \cite{che04b, che04a} becomes the periodic-review inventory control problem,  the technical assumption becomes Assumption~{\bf GB},  and the results in \cite{che04b, che04a} imply the optimality of $(s,S)$ policies.  For coordinating inventory control and pricing, Huh et al.~\cite{huh} provided a method for proving the optimality of stationary policies by adding specific assumptions that hold for inventory control to the MDP assumptions.

   Using the results from Feinberg et al.~\cite{FKZMDP} on the existence of stationary optimal policies and their properties for MDPs with general state and action sets and with possibly unbounded one-step cost functions,
   Feinberg and Lewis~\cite{FL1} proved the optimality of $(s,S)$ policies for a general demand distribution for criteria (i -- iii) mentioned in the previous paragraph.  Feinberg and Liang~\cite{FLi} provided a complete description of optimal discounted policies for arbitrary demand. These results cover the results under Assumption~{\bf GB} as a special case. Feinberg and Liang~\cite{FLi1} proved the validity of the optimality equation for average costs per unit time, while the general results for MDPs \cite{FKZMDP} imply only the validity the optimality inequality.  The conclusions from \cite{FL1, FLi1, FLi} are presented in Section~\ref{s3}.

 Studies of MDPs started with investigations of models with finite state and action sets.  Problems with infinite state and action sets were investigated later.  The two classic objective criteria for infinite-horizon problems are: (i) minimization of  expected total discounted costs, and (ii) minimization of long-run average costs per unit time.  Problems with average cost criteria are usually more difficult.  In particular, optimality equations can be written for expected total costs under mild conditions, and for total expected discounted costs their analyses lead to the proof of optimality of stationary policies for infinite-horizon problems.  For long-run average costs, stationary policies are optimal under stronger conditions than for discounted costs, and proofs of their optimality for average-cost criteria usually use the existence of stationary optimal policies for discounted criteria, when the discount factor increases to 1.  This is the so-called vanishing discount factor approach.  In particular, this approach can be used to establish the validity of optimality equations (sometimes called canonical equations) and inequalities for MDPs with long-run average-costs.  Average-cost optimality equations and inequalities imply the existence of optimal stationary policies for long-run average costs. In applications, average-cost optimality equations and inequalities can be written without an explicit use of the vanishing discount factor approach by using general results on the validity of average-cost optimality equations and inequalities for MDPs. However, as mentioned above, this approach is typically used in the theory of MDPs to establish the validity of such equations and inequalities.

 Let us discuss optimality conditions for  MDPs that are general enough to provide optimality conditions for broad classes  of inventory control models.  First, the state space should be an unbounded subset of a Euclidean space.  This level of generality is covered by Borel state spaces (more precisely, Borel subsets of complete separable metric spaces).  Euclidean spaces are examples of Borel spaces, and the general theory of MDPs with  Euclidean state spaces is not simpler than for Borel spaces. Similarly to subsets of Euclidean spaces, Borel spaces are either finite, countable, or have the cardinality of the continuum. A reader, who is not familiar with the notion of Borel spaces, may view all the state and action sets in this tutorial as subsets of Euclidean spaces. Second, the cost functions may be unbounded.   More precisely, the cost functions should be inf-compact as a function of two variables: a state and action. For inventory control, inf-compact cost functions can be interpreted as lower-semicontinuous functions tending to infinity if either the inventory/backorder or the order size tends to infinity.  Cost functions may not be continuous. For example, they are not continuous in models with positive ordering costs.
      Third, transition probabilities should satisfy the property of continuity in distribution, also known under the name of weak continuity. In particular, transition probabilities are typically weakly continuous for periodic-review stochastic inventory control problems with arbitrary demand distributions;
      see Feinberg and Lewis~\cite[Section 4]{FL} for details.  In particular, it is explained there, that the case of setwise continuous transition probabilities, which is often considered in the MDP literature,  typically covers only discrete and continuous demand distributions.   Fourth, action sets may be unbounded.  This corresponds to a potentially unlimited production/supply capacity.  For example, if a production/supply capacity is limited, then $(s,S)$ policies may not be optimal; see e.g., Federgruen and Zipkin~\cite{FZ} and  Shaoxiang~\cite{Shao}.

 For discounted costs,  Shapley~\cite{Shapley} introduced a zero-sum two-person stochastic game with  finite state and action sets.  If one of the players has only one action at each state, this model becomes an MDP.  This publication is considered as the first paper on MDPs.  Blackwell~\cite{Black} developed the theory for discounted costs and Borel state and action sets.  In particular, Blackwell~\cite{Black} studied problems with bounded costs and discovered that the objective functions may not be Borel measurable, and the dynamic programming approach to such problem should deal with more general policies than Borel measurable ones.  The appropriate theory is developed in Bertsekas and Shreve~\cite{besh96}. Sch\"al~\cite{sch75} developed the theory for discounted costs, Borel state spaces, compact action sets, possibly unbounded above cost functions, and continuous transition probabilities.  Results for two types of continuity are obtained in \cite{sch75}: for setwise and weak continuity. The results on weak continuity are more important for applications and more complicated.  The theory for problems with setwise continuous transition probabilities and possibly noncompact action sets is described in Hern\'{a}ndez-Lerma and Lasserre~\cite{HL:96}.  Feinberg and Lewis~\cite{FL} provided results for discounted MDPs with weakly continuous transition probabilities, possibly uncountable action sets, and inf-compact cost functions. Feinberg et al.~\cite{FKZMDP} introduced the notion of $\K$-inf-compact functions and obtained more general results than in \cite{FL}; see Theorem~\ref{prop:dcoe}, which is a version of \cite[Theorem 2]{FKZMDP} adapted in \cite{FL1} to problems with possibly nonzero terminal costs.  

 For average costs per unit time Blackwell~\cite{Black62} and Derman~\cite{Der} established the existence of stationary optimal policies for the case of finite state and action sets.  Derman~\cite{Der66} and Taylor~\cite{Taylor} introduced optimality equations for infinite-state problems with bounded one-step costs.  These equations and their version for multi-chain problems are called canonical in Dynkin and Yushkevich~\cite{DY}.  Sennott~\cite{sen86} introduced optimality conditions that lead to the validity of optimality inequalities whose solutions define stationary optimal policies; see also \cite{sen99, sen02} and the references therein.  Cavazos-Cadena~\cite{CC} provided an example when optimality inequalities do not hold in the form of equalities. Sch\"al~\cite{sch93} extends Sennott's results to Borel state spaces, compact action spaces, and with weakly and setwise continuous transition probabilities.  Hern\'{a}ndez-Lerma~\cite{hern91} generalized Sch\"al's~\cite{sch93} results for setwise continuous transition probabilities to possibly noncompact action sets. Feinberg and Lewis~\cite{FL} provided sufficient optimality conditions for weakly continuous transition probabilities and possibly noncompact action sets.  Feinberg et al.~\cite{FKZMDP} provided results for weakly continuous transition probabilities that generalize the corresponding results in Sch\"al~\cite{sch93} and Feinberg and Lewis~\cite{FL}; see Subsection~\ref{ss2} below.

 The second topic covered in this tutorial is optimality conditions for POMDPs and, in particular, for inventory control problems with incomplete information on inventory levels.  Research on inventory management with incomplete information was pioneered by Bensoussan et al.~\cite{Bens3, Bens4, BensOx, Bens5}, where particular problems are studied and the existence of optimal policies and convergence of value iterations are established.  In general, for POMDPs there is a well-known reduction, introduced by  Aoki~\cite{Ao}, \AA str\"om~\cite{As},  Dynkin~\cite{Dyn}, and Shiryaev~\cite{sh} of a POMDP to an MDP whose states are posterior probabilities of the states of the original process.   This reduction holds for problems with Borel state, action, and observation sets, and with measurable transition probabilities
 \cite{besh96, HL:96, Rh, Yu}.  However, it provides little information about the existence of optimal policies and the validity of optimality equations.

  This reduction is based on Bayes' formula, which has an explicit
  form only for problems with transition functions that are either
  discrete or have densities. As a result, except the case of finite
  state, action, and observation sets, very little was known on the
  existence of optimal policies for POMDPs. Therefore, the common
  approach is to study applications by problem-specific methods. The
  general approach, applicable to a large variety of applications, for
  verifying optimality conditions for POMDPs  is developed in Feinberg
  et al. \cite{FKZg}, and one of the applications there deals with
  inventory control. The general optimality results on POMDPs are
  presented in Section~\ref{s6}, and an application to inventory
  control is presented in Section~\ref{s7}.

\section{Markov Decision Processes: Definitions and Optimality Conditions}
\label{s2}
An MDP is defined by a tuple $\{\X,\A,\q,c\},$ where   $\X$ is the state space,  $\A$ is the action space,  $\q$ is the transition probability, and $c$ is the one-step cost function.  The state space $\X$ and action space $\A$ are both assumed to be
Borel subsets of Polish (complete separable metric) spaces.
 If an action $a\in \A$ is selected at a state $x\in\X,$ then a cost
$c(x,a)$ is incurred, where $c:\X\times\A\to\overline\R=\R\cup\{+\infty\},$ and the
system moves to the next state according to the probability distribution
$\q(\cdot|x,a)$ on $\X.$  The function $c$ is assumed to be bounded below and
Borel measurable, and $\q$ is a transition probability, that is, $\q(B|x,a)$ is a Borel
function on $\X\times\A$ for each Borel subset $B$ of $\X,$ and $\q(\cdot|x,a)$ is a
probability measure on the Borel $\sigma$-field of $\X$  for each
$(x,a)\in\X\times\A.$

The decision process proceeds as follows: at time $t=0,1,\ldots$ the current
state of the system, $x_\n$, is observed. A decision-maker decides
which action, $a$, to choose, the cost $c(x,a)$ is accrued, the
system moves to the next state  according to $\q(\cdot\given x,a),$
and the process continues. Let $H_\n=(\X\times\A)^{\n}\times\X$ be the
set of histories for $\n=0,1,\ldots\ .$ A (randomized) decision rule at
epoch $\n=0,1,\ldots$ is a regular transition probability $\pi_\n$
from $H_\n$ to $\A.$ 
In other words, (i)
$\pi_\n(\cdot|h_\n)$ is a probability distribution on $\A,$
where $h_\n=(x_0,a_0,x_1,\ldots,a_{\n-1}, x_\n)$
and (ii) for any measurable subset $B \subseteq \A$, the function
$\pi_\n(B|\cdot)$ is measurable on $H_\n.$ A policy $\pi$ is a
sequence $(\pi_0,\pi_1,\ldots)$ of decision rules. Moreover, $\pi$
is called non-randomized if each probability measure
$\pi_\n(\cdot|h_\n)$ is concentrated at one point. A non-randomized
policy is called Markov if all decisions depend only on the current
state and time.  A Markov policy is called stationary if all
decisions depend only on the current state. Thus, a Markov policy
$\phi$ is defined by a sequence $\phi_0,\phi_1,\ldots$ of measurable
mappings $\phi_\n:\X \rightarrow \A.$ 
A stationary policy $\phi$ is defined by a
measurable mapping $\phi:\X \rightarrow \A.$ 

The Ionescu--Tulcea theorem (see \cite[p. 140-141]{besh96} or  \cite[p.
178]{HL:96}) implies that an initial state $x$ and a policy
$\pi$ define a unique probability distribution $\PR_x^\pi$ on the
set of all trajectories $H_\infty=(\X\times \A)^\infty$ endowed with
the product $\sigma$-field defined by the Borel $\sigma$-fields of $\X$
and $\A.$ Let $\E_x^\pi$ be the expectation with respect to this
distribution. For a finite horizon $N=0,1,\ldots$ and a bounded below measurable function
${\bf F}:\X\to{\overline \R}$ called the terminal value, define the
expected total discounted costs
\begin{align}
    v^{\pi}_{N, {\bf F},\alpha}(x) & :=  \E^\pi_x \left[ \sum_{\n=0}^{N-1}
    \alpha^\n c(x_\n,a_\n)+\alpha^N{\bf F}(x_N)\right], \label{discountpi}
\end{align}
where $v_{0,{\bf F},\alpha}^\pi(x)={\bf F}(x),$  $x\in\X,$ $\alpha\ge 0,$  and, if $N=\infty,$ then $\alpha\in [0,1).$ When ${\bf F}(x)=0$ for all $x\in\X,$ we shall write $v^\pi_{N,\alpha}(x)$ instead of
$v^\pi_{N,{\bf F},\alpha}(x).$  When
$N=\infty$ and ${\bf F}(x)=0$ for all $x\in\X$, \eqref{discountpi} defines the infinite horizon expected
total discounted cost of $\pi$ denoted by  $v_\alpha^\pi(x)$ instead of $v_{\infty,\alpha}^\pi(x).$ The
average costs per unit time are defined as
\begin{align}
    w^{\pi}(x) & := \limsup_{N \rightarrow \infty} \frac{1}{N} \E^\pi_x \sum_{\n=0}^{N-1}
    c(x_\n,a_\n). \label{averagepi}
\end{align}
For each function $V^\pi(x)=v_{N,{\bf F},\alpha}^\pi(x)$, $v_{N,\alpha}^\pi(x)$,
$v_\alpha^\pi(x)$, or $w(x)$, define the optimal cost
\begin{align}
V(x) & := \inf_{\pi \in \Pi} V^{\pi}(x), \label{optdisc}
\end{align}
where $\Pi$ is the set of all policies. A policy $\pi$ is called
\textit{optimal} for the respective criterion if $V^\pi(x)=V(x)$ for
all $x\in \X$.

The defined model is too general for the existence of optimal policies.  However, optimal policies exist
under modest conditions, which typically hold for inventory control applications.  The natural conditions
for inventory control applications are that the transition probability $\q$ is weakly continuous and the cost function $c$ is inf-compact.

The transition probability $\q$ is called weakly continuous, if for every bounded continuous function $f:\X\to\R,$ the function
\[\tilde f(x,a):=\int_\X f(y)\q(dy|x,a)           \qquad\qquad x\in\X,\ a\in\A,
\]
is a continuous function on $\X\times\A.$ For an $\overline{\mathbb{R}}$-valued function $f$,  defined on a subset $U$ of a metric space
$\mathbb{U},$ consider the level sets
\begin{equation}\label{def-D}
\mathcal{D}_f(\lambda;U):=\{y\in U \, : \,  f(y)\le
\lambda\},\qquad \lambda\in\R.\end{equation}  A
function $f$ is called \textit{lower semi-continuous} if all the
level sets $\mathcal{D}_f(\lambda;U)$
 are closed, and a function $f$ is called
\textit{inf-compact} 
if all these sets are compact. In particular,  the cost function $c$ is defined on $\U:=\X\times \A$ and the level sets for $c$ are
 \begin{equation}\label{def-Dc}
\mathcal{D}_c(\lambda;\X\times\A)=\{(x,a)\in \X\times\A \, : \,  c(x,a)\le
\lambda\},\qquad \lambda\in\R.\end{equation}

As shown by Feinberg and Lewis~\cite{FL}, for the discounted costs weak continuity of $\q$ and inf-compactness of $c$ imply the existence of optimal policies.  However,  the condition that the function $x$ is inf-compact can be relaxed by considering the class of $\K$-inf-compact functions.

   For two sets $U$ and $V,$  where $U\subset V,$ and for two functions $f$ and $g$ defined on $V$ and $U$ respectively,  function $g$ defined on $U$ is called
 \emph{the restriction of $f$ to $U$} if $g(x)=f(x)$ when $x\in U.$
\begin{defn}{\rm (cp. Definition~\ref{DKKFUN} in Appendix~\ref{Apx}).}\label{DKIXA} Let $\Sp^i$ be metric spaces and $S^i\subseteq \Sp^i,$ $i=1,2.$
A function $f: S^1\times S^2\to \overline{\mathbb{R}}$ is called
$\K$-inf-compact   if, for any nonempty compact subset
 $K$  of $ S^1,$ the restriction of this function to $K\times S^2$ is inf-compact.
\end{defn}

Definition~\ref{DKIXA} corresponds to Definition~\ref{DKKFUN}  of a  $\K$-inf-compact function $u:\Sp^1\times\Sp^2\to \overline{\mathbb{R}}$ on   ${\rm Gr}_{{\Sp^1}}(\Phi)$ in the following way.  For a given function $f: S^1\times S^2\to \overline{\mathbb{R}},$ define $u: \Sp^1\times \Sp^2\to \overline{\mathbb{R}},$
\begin{equation*} u(s^1,s^2):=
\begin{cases} f(s^1,s^2), &\text{if $s^1\in S^1$ and $s^2\in S^2$;}\\
+\infty, &\text{otherwise;}
\end{cases}
\end{equation*}
and $\Phi(s^1):=S^2$ for all $s^1\in \Sp^1.$ Then the function $f: S^1\times S^2\to \overline{\mathbb{R}}$ is $\K$-inf-compact if and only if the function  $u: \Sp^1\times \Sp^2\to \overline{\mathbb{R}}$ is $\K$-inf-compact on ${\rm Gr}_{{\Sp^1}}(\Phi).$  We mainly apply the notion of a $\K$-inf-compact function $f: S^1\times S^2\to \overline{\mathbb{R}}$ to the situation when $S^1=\X,$ $S^2=\A,$ and $\Sp^i,$ $i=1,2,$ are Polish spaces in which the state and action sets $\X$ and $\A$ are defined respectively.  In many inventory control applications, $\Sp^1=\X$ and $\Sp^2=\A.$ So, if the state and action sets $\X$ and $\A$ are explicitly defined as Polish spaces, we assume that  $\Sp^1=\X$ and $\Sp^2=\A$  are Polish spaces containing $\X$ and $\A$  that are mentioned in the definition of an MDP.  The examples include $\X=\R,$ $\X=[0,\infty),$ $\A=\R,$ and $\A= [0,\infty).$

For a function $f:\X\times\A\to \overline{\mathbb{R}},$ $\K$-inf-compactness is a more general and natural property than inf-compactness.  For example, for $\X=\A=\R$ the function $f(x,a)=|x-a|$ is $\K$-inf-compact, but it is not inf-compact.
As shown in Feinberg et al.~\cite{FKZMDP}, the following assumption is sufficient for the existence of optimal policies for discounted MDPs.

\noindent {\bf Assumption~{\bf W*}.} The following conditions hold:

    (i) the transition probability $\q$ is weakly continuous;

    (ii) the cost function $c$ is  $\K$-inf-compact. 

We list some of the properties of MDPs that take place under Assumption~{\bf W*} (see Theorem~\ref{prop:dcoe} for details):
\begin{enumerate}[{1.}]
\item For a  bounded below, lower semi-continuous terminal value function ${\bf F},$ the final-horizon optimality equation holds for all $\a\ge 0:$

 \begin{align}\label{eq43333}
            v_{\n+1,{\bf F},\alpha}(x) & =\min\limits_{a\in \A}\left\{c(x,a)+\alpha
                \int_\X v_{\n,{\bf F},\alpha}(y)\q(dy|x,a)\right\},\quad x\in
                \mathbb{X},\,\,\n=0,1,...,
        \end{align}
        where $v_{0,{\bf F},\alpha}(x)={\bf F}(x)$ for all $x\in \X$.  In particular, this is true for ${\bf F}\equiv 0$ and $v_{0,\alpha}\equiv 0.$

       \item The function $v_\alpha$ is lower semicontinuous, where $\alpha\in [0,1)$.    If the function ${\bf F}$ is bounded below and lower semi-continuous, then the functions $v_{N, {\bf F},\alpha},$ for $N=0,1,\ldots$ and $\alpha\ge 0,$ are  lower semi-continuous. If in addition ${\bf F}(x)\le v_\alpha (x)$ for all $x\in\X,$ then $v_\alpha (x)=\lim_{N\to\infty} v_{N, {\bf F},\alpha} (x),$ where $\alpha\in [0,1)$.  In particular, this is true for ${\bf F}\equiv 0,$ that is, $v_\alpha (x)=\lim_{N\to\infty}v_{N,\alpha} (x),$ where $\alpha\in [0,1)$.

       \item For $\alpha\in [0,1)$ the infinite-horizon value function $v_\alpha$ satisfies the optimality equation
         \begin{align}\label{eq5a}
            v_{\alpha}(x) & =\min\limits_{a\in \A}\left\{c(x,a)+\alpha\int_\X
                v_{\alpha}(y)\q(dy|x,a)\right\},\qquad x\in \X,
        \end{align}
        a stationary optimal policy exists, and a stationary policy $\phi$ is optimal if and only if

        \begin{align}\label{eq51a}
            v_{\alpha}(x) & =c(x,\phi(x))+\alpha\int_\X
                v_{\alpha}(y)\q(dy|x,\phi(x)),\qquad x\in \X.
        \end{align}

      \item If the one-step cost function $c$ is inf-compact, then the value function $v_\alpha$ is inf-compact, when $\alpha\in [0,1).$ The same is true for the value functions $v_{N, {\bf F},\alpha},$ $N=1,2,\ldots,$ when the terminal value ${\bf F}$ is a bounded below, lower semi-continuous function and $\alpha\ge 0.$
\end{enumerate}

In particular, the fourth property is useful for proving the existence of stationary optimal policies for inventory control problems. It is well-known that for average costs per unit time optimal policies may not exist under Assumption~{\bf W*}.  For example,
optimal policies may not exist for a countable state space and finite action sets; see e.g., Ross~\cite[Section 5.1]{Ross} and for a finite state set, compact action sets, and continuous transition probabilities and costs; see e.g., Dynkin and Yushkevich~\cite[Section 7.8]{DY}.  Next we formulate a general condition, that typically holds for inventory control problems, which together with Assumption~{\bf W*} guarantees the existence of optimal policies for average-cost MDPs.  If  $\inf_{x\in
\X}w(x)<+\infty,$  define  for $\alpha\in [0,1)$:
\begin{align*}
    m_{\alpha} & :=\inf\limits_{x\in \X}v_{\alpha}(x),\qquad\qquad
      u_{\alpha}(x):=v_{\alpha}(x)-m_{\alpha}.
\end{align*}

\noindent {\bf Assumption~{\bf B}.} The following conditions hold:

(i) $\inf_{x\in
\X}w(x)<+\infty;$

(ii) $\sup_{\alpha < 1} u_\alpha (x) <\infty$ for all $x\in\X.$

We notice that the function $u_\alpha$ is nonnegative and Assumption~{\bf B} implies that $m_\a$ cannot take infinite values; see Sch\"al~\cite{sch93}. If Assumption~{\bf B(i)}  does not hold then the average-cost problem is trivial: all policies lead to infinite average losses per unit time. This assumption holds in all well-defined problems and usually it is easy to verify.  The validity of Assumption~{\bf B(ii)} probably follows from various ergodicity and communicating conditions, but this relation has not been studied in the literature.  As explained in the text following Theorem~\ref{lemC} below, Assumption~{\bf B(ii)} holds and can be easily verified for inventory control problems.
As shown in Feinberg et al.~\cite{FKZMDP}, Assumptions~{\bf W*} and {\bf B} imply the existence of stationary optimal policies for average-cost MDPs, which follows from the validity of optimality inequalities.

  For $\alpha\in [0,1)$ consider
\begin{align*}
    \underline{w} & =\ilim\limits_{\alpha\uparrow
        1}(1-\alpha)m_{\alpha},\quad\overline{w}=\slim\limits_{\alpha\uparrow
        1}(1-\alpha)m_{\alpha}.
\end{align*}
According to
Sch\"al \cite[Lemma 1.2]{sch93}, Assumption~{\bf B(i)} implies
\begin{align}\label{eq:sch93}
    0\le \underline{w}\le \overline{w}\le w^*< +\infty.
\end{align}
According to Sch\"al \cite[Proposition 1.3]{sch93}, if there exists a measurable
function $u:\ \X\to [0,\infty)$ and a stationary policy $\phi$ satisfying the Optimality Inequality
\begin{align}\label{acoi}
    \underline{w}+u(x) & \geq c(x,\phi(x))+\int u(y)\q(dy|x,\phi(x)), & & x\in
        \X,
\end{align}
then $\phi$ is average-cost optimal and $w(x)=\underline{w}=\overline{w}$ for all $x\in\X$.  Assumptions~{\bf W*} and {\bf B} imply the existence of a stationary policy $\phi$ satisfying optimality inequality  \eqref{acoi}.

Another form of an optimality inequality was introduced in Feinberg et al. \cite{FKZMDP}, where it was shown that, if there exists a measurable
function $u:\X\to [0,+\infty)$ and a stationary policy $\phi$ such
that
\begin{equation}\label{eq7111}
\overline{w}+u(x)\ge c(x, \phi(x))+\int_\X  u(y)\q(dy|x,
\phi(x)),\quad x\in \X,
\end{equation}
then $\phi$ is average-cost optimal and
\begin{equation}\label{eq:7121}
w(x)=w^{\phi}(x)=\slim\limits_{\alpha\uparrow
1}(1-\alpha)v_{\alpha}(x)=\overline{w},\quad x\in \X.
\end{equation}
 Observe that inequality \eqref{eq7111} is weaker than \eqref{acoi} because  \eqref{acoi} implies \eqref{eq7111}.

The existence of stationary optimal policies satisfying inequality \eqref{eq7111} is proved in Feinberg et al.~\cite{FKZMDP} under Assumptions~{\bf W*}  and an assumption called {\bf \underline B} there, which consists of Assumption~{\bf B}(i) and the following assumption  \cite{FKZMDP}:
\begin{equation} \liminf_{\alpha\uparrow 1} u_\alpha(x)<\infty\quad {\rm for\ all}\quad x\in\X, \label{eubbef}
\end{equation}
which is weaker than Assumption~{\bf B(ii)}. However, an example of an MDP, satisfying Assumptions~{\bf W*} and {\bf \underline B} , but 
not satisfying Assumption~{\bf B(ii)}, is currently unknown.

\begin{rem}
{\rm
The definition of an MDP usually includes the sets of available
actions $A(x)\subseteq\A,$ $x\in\X.$  We do not do this explicitly because we allow
$c(x,a)$ to be equal to $+\infty.$  In other words, a feasible pair $(x,a)$ is modeled
as a pair with finite costs.  To transform this model to a one with feasible
action sets, it is sufficient to consider the sets of available actions $A(x)$ such that
$A(x)\supseteq A_c(x),$ where $A_c(x)=\{a\in\A:c(x,a)<+\infty\},$ $x\in\X.$    In
order to transform an MDP with  action sets $A(x)$ to an MDP with the action set $\A,$
it is sufficient to set $c(x,a)=+\infty$ when $a\in\A\setminus A(x),$  $x\in\X.$   Early works on MDPs by
Blackwell~\cite{Black} and Strauch~\cite{Str} considered models with $A(x)=\A$
for all $x\in\X.$ This approach caused some problems with the generality of the
results because the boundedness of the cost function $c$ was assumed and
therefore $c(x,a)\in\R$ for all $(x,a).$  If the cost function is allowed to take
infinitely large values, models with $A(x)=\A$  are as general as models with
$A(x)\subseteq \A,$  $x\in\X.$
}
\end{rem}

\section{MDPs Defined by Stochastic Equations}\label{s3}
Inventory control problems are often defined by equations
\begin{equation}\label{eq:stochs}
x_{\n+1}=F(x_\n,a_\n,D_{\n+1}),\qquad \n=0,1,\ldots,
\end{equation}
where $x_\n$ is the amount of inventory available at the end of  day $t,$ $a_t$ is the ordered quantity at the end of day $t$, and $D_{\n+1}$ is the demand on day $t+1.$ For the classic periodic-review problem with backlogs $F(x,a,D)=x+a-D,$ and for a problem with lost sales $F(x,a,D)=(x+a-D)^+.$ The system can also incur losses of inventory, there could be lead times, and so on.  So, the function $F$ can have a more complicated form, and interpretations of its parameters may be different for different problems.  Also, in this paper we only consider independent and identically distributed demands, that is, $D_1, D_2,\ldots\ .$ are independent and identically distributed.

Let $\Sp$ be a metric space, ${\cal B}(\Sp)$ be its Borel $\sigma$-field, and $\mu$
be a probability measure on   $(\Sp,{\cal B}(\Sp))$.  Consider a stochastic sequence
$x_\n,$ whose dynamics are defined by  equation \eqref{eq:stochs}, where $D_0,D_1,\ldots$ are independent and identically distributed random variables with values in $\Sp$ whose distributions are defined by a probability measure $\mu$ and $F:\X\times\A\times\Sp\to\X$ is a measurable mapping.

  Equation~\eqref{eq:stochs} defines the transition probability
\begin{equation}\label{eq:defqe}
\q(B|x,a)=\int_\Sp 1\{F(x,a,s)\in B\}\mu(ds),\qquad B\in {\cal B}(\Sp),
 \end{equation}
from $\X\times\A\to\X,$  and $\q(\cdot|x_t,a_t) $ is the distribution of $x_{t+1}$  given $x_t$ and $a_t,$ where $1$ is the indicator function.

The following lemma relates Assumption~{\bf W*(ii)} to the problems defined by stochastic equations.
\begin{lem}\label{l:wcq} (Hern\'{a}ndez-Lerma~\cite[p. 92]{HLerma1}). If the function $F$ is continuous then
the transition probability $\q$ is weakly continuous.
\end{lem}
Consider an MDP with the transition probability $\q$ defined by a continuous function $F.$ If the one-step cost function $c$  is inf-compact, then, for a random variable $D$  with the same distribution as $D_1,$ formulae \eqref{eq43333}--\eqref{eq51a} can be rewritten as

 \begin{align}\label{eq43333i}
            v_{\n+1,{\bf F},\alpha}(x) & =\min\limits_{a\in \A}\left\{c(x,a)+\alpha
                \E v_{\n,{\bf F},\alpha}(F(x,a,D))\right\},\quad x\in
                \mathbb{X},\,\,\n=0,1,...,
        \end{align}

         \begin{align}\label{eq5ai}
            v_{\alpha}(x) & =\min\limits_{a\in \A}\left\{c(x,a)+\alpha\E
                v_{\alpha}(F(x,a,D))\right\},\qquad x\in \X,
        \end{align}
        and
        \begin{align}\label{eq51ai}
            v_{\alpha}(x) & =c(x,\phi(x))+\alpha\int_\X
                v_{\alpha}(F(x,\phi(x),D)),\qquad x\in \X.
        \end{align}
Equation \eqref{acoi} becomes

\begin{align}\label{acoii}
    \underline{w}+u(x) & \geq c(x,\phi(x))+\E u(F(x,a,D)), & & x\in
        \X,
\end{align}
and inequality \eqref{eq7111} becomes the same as \eqref{acoii} with \underline{w} replaced with $\overline{w}.$

\section{The Classic Periodic-Review Problem with Backorders}
\label{s4} In this section we consider a
discrete-time periodic-review inventory control problem with back orders and prove
the existence of an optimal $(s,S)$ policy. For this problem the dynamics are
defined by the following stochastic equation
\begin{align}\label{e:invdinv}
x_{\n+1} &=x_\n+a_\n-D_{\n+1}, \quad \n=0,1,2,\ldots,
\end{align}
where $x_\n$ is the inventory at the end of period $\n$, $a_\n$ is the amount  ordered at the end of period $\n,$ and $D_{\n+1}$ is the demand during period $(\n+1)$.  The demand
is assumed to be i.i.d.  In other words, the dynamics of the system is  defined by equation~\eqref{eq:stochs} with the function
$F(x,a,D)=x+a-D.$  Of course, this function is continuous.  Here we consider the case, when there is a single commodity.  In this case, $x_t,$ $a_t,$ and $D_{t+1},$ $t=0,1,\ldots,$ are real numbers.

A decision-maker views the current inventory of a single commodity at the end of the day and
makes an ordering decision. Assuming zero lead times, the products
are immediately available to meet demand. Demand is then realized,
the decision-maker views the remaining inventory, and the process
continues. Assume the unmet demand is backlogged and the cost of
inventory held or backlogged (negative inventory) is modeled as a
convex function. The demand and the order quantity are assumed to be
non-negative. 
The
dynamics of the system are defined by \eqref{e:invdinv}. Let
\begin{enumerate}[(a)]

  \item $\alpha \in (0,1)$ be the discount factor,
  \item $K \geq 0$ be a fixed ordering cost,
  \item $\c> 0$ be the per unit ordering cost,
  \item $D$ be a nonnegative
  random variable with the same distribution as $D_\n,$ and $P(D > 0) > 0,$
  \item $h(\cdot)$ denote the holding/backordering cost per
  period.  It is assumed that $h:\R\to [0,\infty)$ is a convex function, 
   $h(x) \to
  \infty$ as $|x| \to \infty,$ and $\E h(x-D) < \infty$ for all $x \in \R.$
\end{enumerate}
Without loss of generality, assume that $h(0) = 0$. The fact that
$P(D > 0)> 0$ avoids the trivial case.  For example, if $D=0$ almost
surely then the policy that never orders when the inventory level is
non-negative and orders up to zero when the inventory level is
negative, is optimal under the average cost criterion. Note that $\E D < \infty$ since, in view of Jensen's inequality,
 $h(x-\E D)\le\E h(x-D)<\infty.$ 

Let us define the state space $\X=\R,$ the action set $\A=\R^+,$ where $\R^+=[0,\infty),$
the transition probability $\q$ defined in \eqref{eq:defqe} with $F(x,a,D)=x+a-D,$ and
the one-step cost function
\begin{align*}
    c(x,a) & = K 1_{\{a > 0\}} + \c a + \E h(x+a-D).
\end{align*}
The function $c$ is inf-compact and, of course, the function $F$ is continuous.  Therefore, Assumption~{\bf W*} holds.  It is relatively easy to show that Assumption~{\bf B} holds.  Thus, optimality equations exist for finite horizon and infinite horizon problems.  In particular, they exist for problems with total discounted and average-cost criteria.

Optimality equations and inequalities can be written as
%
\begin{align}
v_{\n+1,{\bf F},\alpha}(x)& = \min \{\min_{a > 0} [K + G_{\n,{\bf F},\alpha}(x+a)], 
G_{\n,{\bf F},\alpha}(x)\} - \c x, \label{KDOEF2}\\
v_\alpha(x)& = \min \{\min_{a > 0} [K + G_{\alpha}(x+a)],
G_{\alpha}(x)\} - \c x, \label{KDOE2}\\
w + u(x) & \geq \min \{\min_{a > 0} [K + H(x+a)], H(x)\} - \c x,
\label{KAOE2}
\end{align}
where $t=0,1,\ldots$ and
\begin{align}\label{G111}
G_{\n,{\bf F},\alpha}(x) & := \c x + \E h(x-D) + \alpha \E v_{\n,{\bf F},\alpha}(x-D),\\
    G_{\alpha}(x) & := \c x + \E h(x-D) + \alpha \E v_{\alpha}(x-D),\label{GNG222}\\
    H(x) & := \c x + \E h(x-D) + \E u(x-D).\label{GNG333}
\end{align}
We also write $G_{t,\alpha}$ instead of $G_{\n,{\bf F},\alpha}$ when ${\bf F}\equiv 0.$

\begin{defn}
Let $s_\n$ and $S_\n$ be real numbers such that $s_\n\le S_\n$,
$\n=0,1,\ldots\ .$ Suppose $x_\n$ denotes the current inventory
level at decision epoch $\n$. A policy is called an $(s_\n,S_\n)$
policy at step $\n$ if it orders up to the level $S_\n$ if $x_\n <
s_\n$ and does not order when $x_\n \ge s_\n.$  A Markov 
policy is called an $(s_\n,S_\n)$ policy if it is an $(s_\n,S_\n)$
policy at all steps $\n=0,1,\ldots\ .$ 
A policy is called
an $(s,S)$ policy if it is stationary and it is an $(s,S)$ policy
at all steps $\n=0,1,\ldots\ .$
\end{defn}

The standard methods for proving the optimality of $(s_t,S_t)$ and $(s,S)$ policies for discounted costs was introduced by Scarf~\cite{scarf}, and  is based on the notion of a $K$-convex function.

\begin{defn} 
     A function $f:\R\to\R$ is called $K$-convex, $K \geq 0,$ if for
    each
    $x \leq y$  and for each $\lambda \in (0,1)$,
    \begin{align*}
        f((1- \lambda) x + \lambda y) \leq (1- \lambda) f(x) +
        \lambda f(y) + \lambda K.
    \end{align*}
\end{defn}

For an  inf-compact function $g: \R\to\R,$ let
    \begin{align}
        S & \in \argmin_{x\in\R} \{g(x)\}, \label{defn:S}\\
        s & := \inf \{x\le S\ |\ g(x) \leq K + g(S)\}. \label{defn:s}
    \end{align}
These real numbers exist because the function $g$ is inf-compact.  In addition, $s$ is defined uniquely and does not depend on $S.$  In addition, $s$ is defined uniquely and does not depend on the choice of $S,$ if there are more than one $S$ satisfying \eqref{defn:S}.

The standard method for proving the optimality of $(s_t,S_t)$ policies is to consider $g=G_{N,\alpha},$ $N=1,2,\ldots,$ and prove by induction  that these functions are inf-compact and $K$-convex, which implies from the optimality equation \eqref{KDOEF2} optimality of $(s_t,S_t)$ policies with $S_t$ and $s_t$ defined by \eqref{defn:S}, \eqref{defn:s} with $g=G_{t,\alpha}$.  The next step would be to consider $t\to\infty$ and  prove the optimality of $(s,S)$ policies for infinite-horizon problems.

However, it is possible that the functions  $G_{N,\alpha}$ are not inf-compact, and the described approach fails. Then the natural approach is to try to do the same steps for the function  $G_{N,{\bf F},\alpha}$ for a specially selected terminal value function ${\bf F}.$  The natural candidate is  the function ${\bf F}=v^0_\alpha,$ where $v^0_\alpha$ is the infinite-horizon value for the problem with the ordering cost $K=0.$ It is possible to show that there exists $\alpha^\prime\in [0,1)$ such that the functions  $G_{N,{v^0_\alpha},\alpha}$ are inf-compact, and this implies the optimality of $(s_t,S_t)$-policies for all finite-horizon problems with the terminal value ${\bf F}=v^0_\alpha$ for all $\alpha\in [\alpha^\prime,1),$ which implies optimality of $(s,S)$-policies for the infinite horizon discounted criterion with the discount factor $\alpha.$ In addition, it is always true that $G_{N,{v^0_\alpha},\alpha}\to G_\alpha$, and the following lemma holds.
 \begin{lem}\label{prop:finite-exp} (\cite{vei65, FL1}).
    There exists $\alpha^\prime\in [0,1)$ such that $G_\alpha(x)\to\infty$ as $|x|\to\infty$ for all $\alpha\in[\alpha^\prime,1) $ and for all setup costs $K\ge 0.$
\end{lem}

 The optimality of $(s,S)$-optimal policies for large discount factors imply optimality of $(s,S)$ policies for average costs per unit time.  The following theorem takes place.

\begin{thm} \label{th:sSdisc} (\cite{FL1}).
Consider $\alpha^\prime\in [0,1)$  whose existence is stated in
Lemma~\ref{prop:finite-exp}. The following statements hold for the inventory
control problem.

   (i) For $\alpha \in [\alpha^\prime,1)$ and $\n=0,1,\ldots,$ define
          $g(x):=G_{\n,v_\alpha^0,\alpha}(x),$ $x\in\R$. Consider real numbers
          $S^*_{\n,\alpha}$ satisfying  \eqref{defn:S} and $s^*_{\n,\alpha}$
          defined in   \eqref{defn:s}. Then for each $N=1,2,\ldots,$  the
          $(s^*_{N-\n,\alpha},S^*_{N-\n,\alpha})$ policy, $\n=1,2,\ldots,N,$ is
          optimal for the $N$-horizon problem with the terminal values
          ${\bf F}(x)=v^0_\alpha(x),$ $x\in\R$. \label{state:inv1}

 (ii)      For the infinite-horizon expected total discounted cost criterion with a
          discount factor $\alpha\in [\alpha^\prime,1),$ define $g(x):=G_\alpha(x),$
          $x\in\R$. Consider real numbers $S_{\alpha}$ satisfying \eqref{defn:S}
          and $s_{\alpha}$ defined in   \eqref{defn:s}. Then the
          $(s_\alpha,S_\alpha)$ policy is optimal for the discount factor $\alpha.$
          Furthermore, the sequence of pairs
          $\{(s^*_{\n,\alpha},S^*_{\n,\alpha})\}_{\n=0,1,\ldots} $ is bounded, where $s^*_{\n,\alpha}$ and $S^*_{\n,\alpha}$
          are described in statement (i), $\n=0,1,\ldots\ .$ If
          $(s^*_\alpha,S^*_\alpha)$  a limit point  of this sequence, then the
          $(s^*_\alpha,S^*_\alpha)$ policy is optimal for the infinite-horizon problem with the discount factor
          $\alpha.$ \label{state:inv2}

(iii) Consider the infinite-horizon average cost  criterion. For each $\alpha\in
          [\alpha^\prime,1)$, consider an optimal $(s^\prime_\alpha, S^\prime_\alpha)$
          policy for the discounted cost criterion with the discount factor $\alpha,$
          whose existence follows from Statement (ii).   Let $\alpha_\n\uparrow 1,$
          $\n=1,2,\ldots,$ with $\alpha_1\ge\alpha^\prime.$  Every sequence
          $\{(s^\prime_{\alpha_\n}, S^\prime_{\alpha_\n})\}_{\n=1,2,\ldots}$ is
          bounded and  each limit point $(s^\prime,S^\prime)$ defines an
          average-cost optimal $(s^\prime,S^\prime)$ policy. \label{state:inv3}
\end{thm}

As explained above, $(s_t,S_t)$ policies may not be optimal for finite-horizon problems for all discount factors and $(s,S)$
may not be optimal for infinite-horizon discounted problems with a small discount factor.  Let us consider the assumption on the growth of backordering costs, that was probably introduced by  Veinott and Wagner~\cite{vei65} for problems with discrete demand.  This assumption ensures that the functions $G_{N,\alpha}$ and $G_\alpha$ are inf-compact, and, as explained above, this implies the optimality of $(s_t,S_t)$ policies and $(s,S)$ policies for finite-horizon and infinite-horizon discounted problems respectively for all $N=1,2,\ldots$ and for all $\alpha\in [0,1).$

\noindent{\bf Assumption~{\bf GB}.}
There exist $z,y\in\R$ such that $z<y$ and
\begin{align}
                \frac{\E [h(y-D) - h(z-D)]}{y-z} < -\c. \label{eq:slope}
            \end{align}

\begin{lem}\label{lem:galpha} (\cite{che04b, che04a, FL1, heysol}).
       Suppose that Assumption~{\bf GB} holds.
        Then the functions $G_\alpha(x)$  and $G_{N,\alpha}(x),$ $N=1,2,\ldots,$ are inf-compact and $K$-convex.
\end{lem}
The following theorem describes the optimality of $(s_t,S_t)$ policies and $(s,S)$ policies for finite-horizon and infinite-horizon discounted problems under Assumption~{\bf GB}.
\begin{thm}\label{th:sSdisc1} (\cite{che04b, che04a, FL1}). Suppose that  Assumption~{\bf GB} holds.  Then:

(i) For $\alpha \ge 0$ and $\n=0,1,\ldots,$ consider real numbers
      $S_{\n,\alpha}$ satisfying \eqref{defn:S} and $s_{\n,\alpha}$ defined in
      \eqref{defn:s}  with $g(x)=G_{\n,\alpha}(x),$ $x\in\R.$ Then for
      every $N=1,2,\ldots$  the $(s_{N-\n,\alpha},S_{N-\n,\alpha})$ policy,
      $\n=1,2,\ldots,N,$ is an optimal policy for the $N$-horizon problem with the
      zero terminal values.

 (ii) Let $\alpha\in [0,1).$ Consider real numbers $S_{\alpha}$ satisfying \eqref{defn:S}
          and $s_{\alpha}$ defined in   \eqref{defn:s} for $g(x):=G_\alpha(x),$
          $x\in\R.$ Then the
          $(s_\alpha,S_\alpha)$ policy is optimal for the infinite-horizon problem with the discount factor $\alpha.$ 
 Furthermore, a sequence of pairs
          $\{(s_{\n,\alpha},S_{\n,\alpha})\}_{\n=0,1,\ldots} $ considered in statement (i) is bounded, and, if
          $(s^*_\alpha,S^*_\alpha)$  is a limit point of this sequence, then the
          $(s^*_\alpha,S^*_\alpha)$ policy is optimal for the infinite-horizon problem with the discount factor
          $\alpha.$
\end{thm}

As stated in Theorem~\ref{th:sSdisc}, $(s,S)$-policies are optimal for average costs per unit time.  However, Theorem~\ref{th:sSdisc1} states the
optimality of $(s_t,S_t)$ policies and $(s,S)$ policies for finite-horizon and infinite-horizon discounted problems for all discount factors  only under Assumption~{\bf GB}.  The structure of discount optimal policies for all discount factors is investigated in Feinberg and Liang~\cite{FLi}, where the following parameters were introduced:
\begin{align}
	k_h:= - \lim_{x\to -\infty} \frac{h(x)}{x} .
	\label{eqn:limit h(x)}
\end{align}
and
\begin{align}
	\a^* := 1 - \frac{k_h}{\c} .
	\label{def a*}
\end{align}
For example, $\a^*=1-\frac{h_-}{\c}$ for models with linear holding and bacordering costs $h$ considered in \cite{Bens0, bert00}, when
\[
h(x)=\begin{cases}\ \ h_+x, &\text{if\ $x\ge 0;$}\\ -h_- x,& {\rm otherwise;}
\end{cases}
\]
where $h_-$ and $h_+$ are positive holding and backordering cost rates, and typically $h_->h_+.$

The convexity and inf-compactness of $h$ imply that  $0< k_h \leq +\infty$.  Therefore,  $-\infty\leq \a^* < 1$. In addition, Assumption~{\bf GB} is equivalent to $\a^* < 0$. In addition, $\a^\prime:=\max\{\a^*,0\}$ is the minimal possible value of the parameter $\alpha^\prime$ whose existence is claimed in Lemma~\ref{prop:finite-exp}.  These facts and their corollaries are summarized in the following theorem.

\begin{thm} (\cite{FLi}). Assumption~{\bf GB} holds if and only if $\alpha^*<0.$ Therefore, if $\alpha^*<0,$ then the statements (i) and (ii) of Theorem~\ref{th:sSdisc1} hold.  In addition, $\a^\prime=\max\{\a^*,0\}$ is the minimal value of the parameter $\a^\prime$ whose existence is stated in Lemma~\ref{prop:finite-exp}. Therefore, statements (i) and (ii) of Theorem~\ref{th:sSdisc} take place for $\alpha^\prime=\max\{\a^*,0\}.$
\end{thm}

Define  $\bS_0:=0$ and
\begin{align}
	\bS_t := \sum_{j=1}^t D_j,\qquad\qquad\qquad t=1,2,\ldots\ .
	\label{eqn:sum of demand}
\end{align}
 Then $\E[\bS_t] = t\E[D] < +\infty$
for all $t=0,1,\ldots\ .$

Define the following function for all $t=0,1,\ldots$ and $\a\geq 0,$
\begin{align}
	f_{t,\a}(x) := \c x + \sum_{i=0}^{t} \a^i \E[h(x-\bS_{i+1})] , \quad x\in\X.
	\label{eqn:fta(x)}
\end{align}
Observe that $f_{0,\a} (x)= \c x + \E[h(x-D)] = G_{0,\a} .$ Since $h(x)$ is a convex function,
then the function $f_{t,\a}(x)$ is convex for all $t=0,1,\ldots$ and $\a\geq 0.$

Let $F_{t,\a}(-\infty):=\lim_{x\to -\infty}f_{t,\a}(x)$ and
\begin{align}
	N_{\a} := \inf \{ t=0,1,\ldots: F_{t,\a}(-\infty) = +\infty \},
	\label{eqn:def Na}
\end{align}
where the infimum of an empty set is $+\infty$. Since the function $h(x)$ is non-negative,
then the function $f_{t,\a}(x)$ is non-decreasing in $t$ for all $x\in\X$ and $\a\geq 0.$
Therefore, (i) $N_{\a}$ is non-increasing in $\a$, that is, $N_{\a}\leq N_{\beta}$,
if $\a>\beta;$ and (ii) in view of the definition of $N_{\a}$, for each $t\in\N_0$
\begin{align}
	F_{t,\a}(-\infty)  < +\infty , \qquad \text{if } t < N_{\a} , \qquad \text{and } \qquad
	F_{t,\a}(-\infty)	= +\infty, \qquad \text{if } t \geq N_{\a} .
	\label{eqn:limit of fta}
\end{align}

The following theorem provides the complete description of optimal finite-horizon  policies for all discount factors $\alpha.$
\begin{thm}\label{thm:general results}  (\cite{FLi}).  Let $\alpha>0.$
	Consider $\a^*$ defined in \eqref{def a*}. If $\a^*<0$ (that is, Assumption~{\bf GB} holds), then the statement of Theorem \ref{th:sSdisc1}(i)
	holds. If $0\leq \a^* < 1,$ then the following statements hold for the finite-horizon problem with the discount factor $\alpha:$
	\begin{enumerate}[(i)]
		\item if $\a\in [0,\a^*],$ then a policy that never orders is optimal for
		every finite horizon $N=1,2,\ldots;$
\item if $\a > \a^*$, then $N_{\a}<+\infty$ and for a finite horizon
		$N=1,2,\ldots,$  the following is true:
	\begin{enumerate}[(a)]
			\item if $N \leq N_{\a},$ then
			a policy that never orders at steps $t = 0,1,\ldots,N-1$ is optimal;
			
\item if $N> N_{\a},$ then a policy that never orders at steps
			$t = N-N_{\a},\ldots,N-1$ and follows the $(s_{N-t-1,\a},S_{N-t-1,\a})$ policy
			at steps $t = 0,\ldots,N-N_{\a}-1$ is optimal, where the real numbers $S_{t,\a}$
			satisfy \eqref{defn:S} and  	$s_{t,\a}$ are defined in
			\eqref{defn:s}  with $g(x):=G_{t,\a}(x),$ $x\in\X$.
\end{enumerate}\end{enumerate}
\end{thm}

\begin{figure}[h]
  \centering
  \caption{The structure of optimal policies for a discounted $N$-horizon problem with $N < +\infty$ and $\alpha \geq 0$.}
  \label{fig:optPol_finite}
  \includegraphics[scale=0.3]{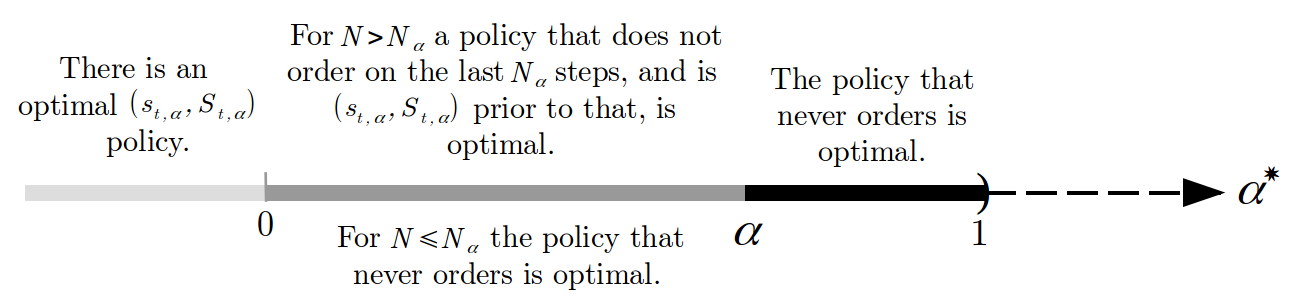}
\end{figure}

The conclusions of Theorem~\ref{thm:general results}  are presented in Table~\ref{Table1}   and Figure~\ref{fig:optPol_finite}.

\begin{table}[ht]
\centering
\caption{The structure of optimal policies for a discounted $N$-horizon problem with $N<+\infty$ and $\alpha\ge 0.$}\label{Table1}

	\begin{tabular}{|c|l|l|l|}
    	\hline
    $\alpha$ & \multicolumn{1}{|c|} {$\alpha^* < 0$} &  \multicolumn{1}{|c|}{$0\leq \alpha^*<\alpha$} & \multicolumn{1}{|c|}  {$\alpha^* \geq \alpha$} \\
	\hline
& There is  & For the natural number $N_{\a}$ defined in \eqref{eqn:def Na},  & The policy\\
& an optimal & $\ \ \ $ if $N > N_{\alpha},$ then a policy that never orders at & that never \\
& $(s_{t,\a},S_{t,\a})$   & steps $t = N-N_{\alpha},\ldots,N-1$ and is an $(s_{t,\a},S_{t,\a})$ & orders is  \\
& policy. & policy at steps $t = 0,\ldots,N-N_{\alpha}-1$ is optimal; &  optimal. \\
&	&  $\ \ \ $ if $N \leq N_{\alpha},$ then a policy that never orders is  & \\
&	&  optimal. & \\
		\hline
	\end{tabular}
	\label{tab:optimal policies T horizon}

\end{table}

The following theorem provides the complete description of optimal infinite-horizon  policies for all discount factors $\alpha.$

\begin{thm}\label{thm:general results_ig} (\cite{FLi}). Let $\alpha\in [0,1).$
	Consider $\a^*$ defined in \eqref{def a*}. The following statements hold for the infinite-horizon problem with the discount factor $\alpha:$

(i) if $\a^*<\a,$ then an $(s_\alpha, S_\alpha)$ policy is optimal, where the real numbers $S_{\a}$ satisfy
			 \eqref{defn:S} and 	$s_{a}$ are defined in
			\eqref{defn:s}  with $g(x):=G_{\a}(x),$ $x\in\X.$ Furthermore, a sequence of pairs ${(s_{t,\a},S_{t,\a})}_{t=N_{\a},N_{\a}+1,\ldots}$ considered in Theorem \ref{thm:general results} (ii,b) is bounded, and, for if
$(s_{\a}^*,S_{\a}^*)$ is a limit point of the sequence, then the $(s_{\a}^*,S_{\a}^*)$ policy is optimal for the infinite-horizon problem with the discount factor $\a;$

(ii) if $\alpha^*\ge\alpha,$ then the policy that never orders is optimal.
\end{thm}

\begin{figure}[h]
  \centering
  \caption{The structure of optimal policies for a discounted infinite-horizon problem with $\alpha \in [0, 1)$.}
  \label{fig:optPol_infinite}
  \includegraphics[scale=0.3]{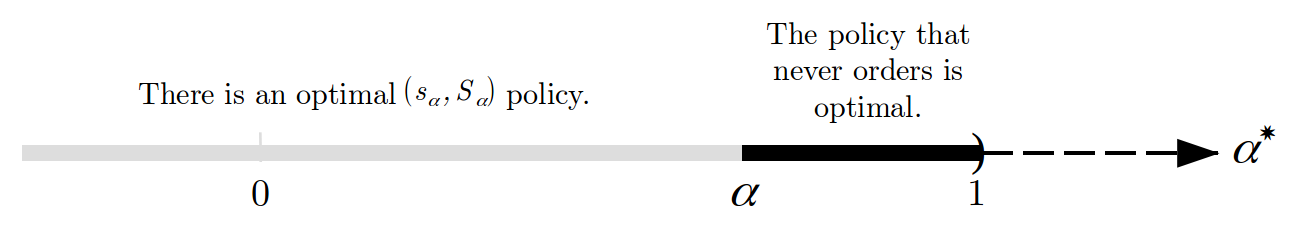}
\end{figure}

The conclusions of Theorem~\ref{thm:general results_ig}  are presented in Table~\ref{Table2} and Figure~\ref{fig:optPol_infinite}.

\begin{table}[ht]
\centering
\caption{The structure of optimal policies for a discounted infinite-horizon problem with $\alpha\in [0,1)$.}\label{Table2}
	\begin{tabular}{|c|l|l|}
    	\hline
    $\a$ & \multicolumn{1}{|c|}{$\a^*<\a$} & \multicolumn{1}{|c|}{$\alpha \leq \alpha^*$} \\
	\hline
	& There is an optimal  & The policy that never \\
	& $(s_{\a},S_{\a})$ policy. &  orders is optimal. \\
	\hline
	\end{tabular}
	\label{tab:optimal policies infinite horizon}
	
\end{table}

The above theorems describe stationary optimal policies for all discount factors.  However, it is possible that for a given discount factor at some states there are multiple optimal actions.  Therefore, there may exist multiple stationary optimal policies. It is also possible to describe all stationary optimal policies; Feinberg and Liang~\cite{FLi}.  The results on MDPs imply that the functions $v_\alpha$ and $v_{N,\alpha}$ are lower semi-continuous.  However, for this problem they are continuous; Feinberg and Liang~\cite{FLi}.  In addition, optimality inequalities \eqref{acoi}  and \eqref{KAOE2} hold in the form of equalities; Feinberg and Liang~\cite{FLi1}.

\section{MDPs with Infinite State Spaces and Weakly Continuous Transition Probabilities}
This section describes the theory of dynamic programming for infinite-state problems with  weakly continuous transition probabilities.  The main focus is on the existence of optimal policies and the validity of optimality equations for problems with discounted costs and optimality inequalities for average-cost problems.  We also discuss the convergence of optimal values and actions when the horizon length tends to infinity for finite horizon problems and when the discount factor increases to 1 for infinite horizon problems.

\subsection{Total Discounted Costs}
The following theorem describes the validity of optimality equalities, the lower semi-continuity of value functions and the convergence of value iterations.  For zero terminal values, this theorem is presented in Feinberg et al.~\cite{FKZMDP}.  The case of nonzero terminal values is added in Feinberg and Lewis~\cite{FL1}.  The case of inf-compact cost functions $c,$ which leads to the inf-compactess of value functions, is studied in Feinberg and Lewis~\cite{FL}. The inf-compactness of value functions is important for the analysis of average-cost problems.  The proof of Theorem~\ref{prop:dcoe} uses the generalization of Berge's theorem described in Appendix.

\begin{thm}\label{prop:dcoe} (\cite{FKZMDP, FL1}).
    Let Assumption~{\bf W*} hold.  Consider a  bounded below, lower
    semi-continuous function ${\bf F}:\X\to\overline\R$ and $\a\ge 0.$ Then:
\begin{enumerate}[(i)]
    \item the functions $v_{\n,{\bf F},\alpha}$, $\n=0,1,\ldots,$ are
          lower semi-continuous; \label{cor:finite1}

        \item  the finite-horizon optimality equalities \eqref{eq43333} hold with $v_{0,{\bf F},\alpha}(x)={\bf F}(x)$ for all $x\in \X$
         and the nonempty sets
        \begin{align*}
            A_{\n,{\bf F},\alpha}(x) & :=\{a\in \A:\,v_{\n+1,{\bf F},\alpha}(x)=c(x,a)+\alpha \int_\X
                v_{\n,{\bf F},\alpha}(y)\q(dy|x,a) \},\ \ x\in \X, \n=0,1,\ldots,
        \end{align*}
        satisfy the following properties:
\begin{enumerate}[(a)]
       \item the graph ${{\Gr}}_\X(A_{\n,{\bf F},\alpha})=\{(x,a):\, x\in\X, a\in
              A_{\n,{\bf F},\alpha}(x)\}$, $\n=0,1,\ldots,$ is a Borel subset of $\X\times
              \mathbb{A}$, and

        \item if $v_{\n+1,{\bf F},\alpha}(x)=+\infty$, then $A_{\n,{\bf F},\alpha}(x)=\A$ and, if
              $v_{\n+1,{\bf F},\alpha}(x)<+\infty$, then $A_{\n,{\bf F},\alpha}(x)$ is compact;
\end{enumerate}
       \item for a problem with the terminal value function ${\bf F},$ for each
          $N=1,2,\ldots$, there exists a Markov optimal $N$-horizon policy
          $(\phi_0,\ldots,\phi_{N-1})$ and if, for an $N$-horizon Markov policy
          $(\phi_0,\ldots,\phi_{N-1})$ the inclusions $\phi_{N-1-\n}(x)\in
          A_{\n,{\bf F},\alpha}(x)$, $x\in\X,$ $\n=0,\ldots,N-1,$ hold then this policy is
          $N$-horizon optimal; 

 \item if the cost function $c$ is inf-compact, the functions $v_{\n,{\bf F},\alpha},$
          $\n=1,2,\ldots,$ are inf-compact. 

   \item  for $\a\in [0,1),$ if ${\bf F}(x)$ is constant or ${\bf F}(x)\le v_\alpha(x)$ for all $x\in\X,$ then $v_{\n,{\bf F},\alpha}(x)\to
          v_\alpha (x)$ as $\n \to +\infty$ for all $x\in \X;$ 

     \item       for $\alpha\in [0,1),$ the infinite-horizon optimality equation
        \eqref {eq5a}
           holds
        and the nonempty sets \[A_{\alpha}(x):=\{a\in
        \A:\,v_{\alpha}(x)=c(x,a)+\alpha\int_\X
                v_{\alpha}(y)\q(dy|x,a) \},\qquad x\in \X,\] satisfy
        the following properties: 
\begin{enumerate}[(a)]
       \item    the graph ${\Gr}_\X(A_{\alpha})=\{(x,a):\, x\in\X, a\in
                A_\alpha(x)\}$ is a Borel subset of $\X\times \mathbb{A}$, and

 \item if $v_{\alpha}(x)=+\infty$, then $A_{\alpha}(x)=\A$ and, if
                $v_{\alpha}(x)<+\infty$, then $A_{\alpha}(x)$ is compact.
\end{enumerate}
   \item for an infinite-horizon problem with $\alpha\in [0,1)$ there exists a stationary discount-optimal policy
          $\phi_\alpha$, and a stationary policy $\phi_\a$ is optimal if and only if
          $\phi_\alpha(x)\in A_\alpha(x)$ for all $x\in \X.$

 \item if the cost
          function $c$ is inf-compact, then  the infinite-horizon value function $v_\alpha$ is inf-compact, $\alpha\in [0,1).$ \label{state:fl1}
\end{enumerate}
\end{thm}

The following theorem describes convergence properties of  optimal finite-horizon actions as the time horizon increases to infinity.
\begin{thm}\label{t:lima}  (\cite{FL1}.) Let Assumption~{\bf W*} hold and $\alpha\in [0,1).$ Let ${\bf F}:\X\to\overline \R$ be bounded below, lower semi-continuous, and such that  for all $x\in\X$
\begin{equation}\label{eq:spdouco}
 {\bf F}(x)\le v_\alpha(x)\qquad{\rm and} \qquad v_{1,{\bf F},\alpha}(x)\ge {\bf F}(x).
\end{equation}
Then for $x\in\X,$ such that $v_\alpha(x)<\infty,$ the following two statements hold:
  \begin{enumerate}[(i)]
\item
 there is a compact subset $D^*_\alpha(x)$ of $\A$ such that $A_{\n,{\bf F},\alpha}(x)\subseteq D^*_\alpha(x)$ for all $t=1,2,\ldots,$  where the sets $A_{\n,{\bf F},\alpha}(x)$ are defined in  Theorem~\ref{prop:dcoe}(ii);\item
  each sequence $\{a^{(\n)}\in A_{\n,{\bf F},\alpha}(x)\}_{\n=1,2,\ldots}$  is bounded, and all its limit points belong to $A_\alpha(x).$ \end{enumerate}
\end{thm}

Theorem~\ref{t:lima} is useful for the analysis of the classic periodic-review inventory problem described in Section~\ref{s4}.  As demonstrated in Table~\ref{Table1}, $(s_t,S_t)$ policies may not be optimal for finite horizon problems, and   the function ${\bf F}=v_{0}^0$ is used   to approximate optimal infinite-horizon thresholds, where $v^0_\alpha$ is the optimal value in the same problem with zero ordering costs.
\subsection{Average Costs per Unit Time} \label{ss2}
We start with the formal introduction of Assumption~{\bf\underline B}.

\noindent {\bf Assumption~{\bf \underline B}.} The following conditions hold:

(i) $\inf_{x\in
\X}w(x)<+\infty;$

(ii) $\liminf_{\alpha < 1} u_\alpha (x) <\infty$ for all $x\in\X.$

Recall that the functions $v_\alpha$ and $u_\alpha$ are defined only for $\alpha\in [0,1).$ Let us set
\begin{equation}\label{eq131}
u(x):=\ilim\limits_{ (y,\alpha)\to (x,1-)}u_{\alpha}(y),\qquad\qquad  x\in\X.
\end{equation}
In words, ${u}(x)$ is the largest number such that ${u}(x)\le \liminf_{n\to\infty}u_{\alpha_n}(y_n)$ for all sequences $\{y_n\to x\}$ and $\{\alpha_n\to 1-\}.$

\begin{thm}\label{teor1} (Feinberg et al. \cite[Theorem 3]{FKZMDP}).
Suppose Assumptions~{\bf W*} and {\bf
\underline B} hold. Then there exists a stationary policy $\phi$
satisfying (\ref{eq7111}) with $u$ defined in (\ref{eq131}). Thus,
equalities (\ref{eq:7121}) hold for this policy $\phi.$
Furthermore, the following statements hold:
\begin{enumerate}[(i)]
\item the function $u:\X\to \R_+$ is lower semi-continuous;

\item the nonempty sets
\begin{equation}\label{defsetA*u}
 A_u^*(x):=\left\{a\in \A\, : \,
\overline{w}+u(x)\ge c(x,a)+\int_\X  u(y)\q(dy|x,a) \right\}, \
x\in\X,
\end{equation}
satisfy the following properties:
\begin{enumerate}[(a)]
\item the graph ${\rm Gr}(A_u^*)=\{(x,a):\, x\in\X, a\in
A_u^*(x)\}$ is a Borel subset of $\X\times \mathbb{A}$;

\item  for each $x\in\X$ the set $A_u^*(x)$ is compact;
\end{enumerate}
\item a stationary policy $\phi$ is optimal for
average costs and satisfies (\ref{eq7111}) with $u$ defined in
(\ref{eq131}), if $\phi(x)\in A_u^*(x)$ for all $x\in \X$;
 \item  there exists a stationary policy $\phi$ with
 $\phi(x)\in A_*(x)\subseteq
A_u^*(x)$ for all $x\in\X$, where
\begin{equation}\label{defsetA*}
 A_*(x):=\left\{a\in \A\, : \,
u(x)=\inf_{a\in\A}\{ c(x,a)+\int_\X  u(y)\q(dy|x,a)\} \right\}, \qquad
x\in\X,
\end{equation}
 \item
if, in addition, the function $c$ is inf-compact, then the
function $u$  is inf-compact. 
\end{enumerate}
\end{thm}

Stronger results hold under Assumption~{\bf B}.

\begin{thm}\label{teor3} {\rm (Feinberg et al.~\cite[Theorem 4]{FKZMDP})}.
    Suppose Assumptions~{\bf W*} and {\bf B} hold. Then there exists a nonnegative lower
    semi-continuous function $u$ and a stationary policy $\phi$ satisfying \eqref{acoi},
    that is, $\phi(x)\in A^*_u(x)$ for all $x\in\X$. Furthermore, every stationary policy
    $\phi$, for which \eqref{acoi} holds, is optimal for the average costs per unit time
    criterion,
    \begin{align} \label{eq5.16}
        w^\phi(x) & =w(x)=w^*=\underline{w}=\overline{w}=\lim\limits_{\alpha\uparrow
            1}(1-\alpha)v_\alpha(x)=\lim\limits_{N\to\infty}\frac{1}{N}v^\phi_{N,1}(x),\qquad
            x\in\X.
    \end{align}
    Moreover, the following statements hold:
    \begin{enumerate}[(i)]
     \item the nonempty sets $A_u^*(x), x\in\X$, satisfy the following properties:
     \begin{enumerate}[(a)]
   \item the graph ${{\Gr}_\X}(A_u^*)=\{(x,a):\, x\in\X, a\in A_u^*(x)\}$ is a
              Borel subset of $\X\times \mathbb{A}$;
  \item        for each $x\in\X$ the set $A_u^*(x)$ is compact;
          \end{enumerate}
\item there exists a stationary policy $\phi$ with $\phi(x)\in A_u^*(x)$ for all
          $x\in\X.$
    \end{enumerate}
\end{thm}

Alternatively to \eqref{eq131}, as follows from Feinberg et al.~\cite[Theorems 3,4 and p. 603]{FKZMDP},  for each sequence $\alpha_n\to 1-,$ the function $u$ can be defined as
\begin{align}\label{eq1311}
    {\tilde u}(x) & :=\ilim\limits_{ (y,n)\to (x,\infty)}u_{\alpha_n}(y),\quad x\in\X.
\end{align}
In words, ${\tilde u}(x)$ is the largest number such that ${\tilde u}(x)\le \liminf_{n\to\infty}u_{\alpha_n}(y_n)$ for all sequences $\{y_n\to x\}.$ 
It follows from these definitions that $u(x)\le{\tilde u}(x),$ $x\in\X.$ However, the
questions, whether $u={\tilde u}$ and whether the values of $\tilde u$ depend on a particular choice of the sequence $\alpha_n$ has not been investigated. If the cost function $c$
is inf-compact, then the functions $v_\alpha,$  $u,$ and $\tilde u$ are inf-compact
as well; see Theorem~\ref{prop:dcoe} for the proof of this fact for $v_\alpha$ and
Feinberg et al.~\cite[Theorem 4(e) and Corollary 2]{FKZMDP} for $u$ and $\tilde
u$. We denote by $A^*_{\tilde u}(x)$ the sets  defined
          in \eqref{eq131}, when the function $u$ is replaced with $\tilde u.$

In addition, if the one-step cost function $c$ is inf-compact,  the minima of the
functions $v_\alpha$ possess additional properties.
Set
\begin{align}\label{e:defxa}
    X_\alpha & :=\{x\in\X\,: v_\alpha(x)=m_\alpha\},\qquad \alpha\in [0,1).
\end{align}
In view of Theorem~\ref{prop:dcoe}(viii), the function $v_\a$ is inf-compact and $X_\alpha\ne\emptyset.$  Since $X_\alpha=\{x\in\X\,: v_\alpha(x)\le m_\alpha\},$ this set is closed.    The following fact is useful for
verifying the validity of Assumption~{\bf B(ii)} in inventory control applications; see Feinberg and Lewis~\cite[Lemma
5.1]{FL} and the references therein.

\begin{thm}\label{lemC} (Feinberg et al.~\cite[Theorem 6]{FKZMDP}).
Let Assumptions~{\bf W*} and {\bf B(i)} hold. If the function $c$ is inf-compact, then
there exists a compact set $\mathcal{K}\subseteq\X$ such that
$X_\alpha\subseteq \mathcal{K}$ for all $\alpha\in [0,1).$
\end{thm}

\begin{figure}[h]
  \centering
  \caption{Verification of Assumption~{\bf B(ii)} for inventory control.}
  \label{verify_asmpB}
  \includegraphics[scale=0.3]{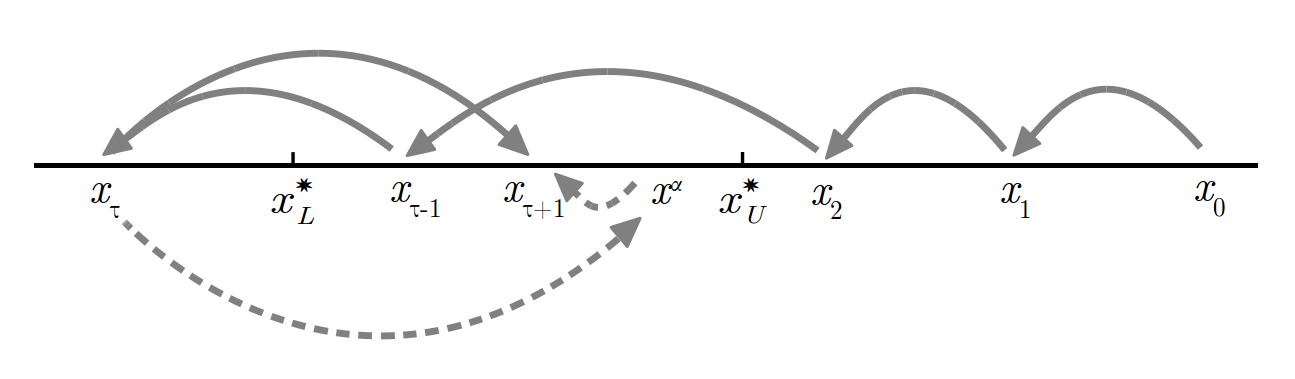}
\end{figure}

Theorem~\eqref{lemC} implies that the minimum in $x\in\X$ of $v_\alpha(x)$ is achieved on a compact set $\mathcal{K},$ which does not depend on $\alpha.$
This typically means that to prove Assumption~{\bf B(ii)} it is sufficient to show that for each $x\in\X$ it is possible to reach every point in $\mathcal{K}$ in a way that the expected time and cost are finite.  In inventory control applications this can be shown  by lowering the inventory levels below the levels in $\mathcal{K}$ and then by ordering up to a point in $\mathcal{K}.$  Exact mathematical justifications are usually problem-specific and use  renewal theory.  Here we provide a short version of the proof from Feinberg and Lewis~\cite{FL1}.  Choose  ${\cal K} =[x^*_L,x^*_U];$ see Figure~\ref{verify_asmpB}, where the existence of a set $\cal K$ is stated in Theorem~\ref{lemC}, and this set can be chosen to be equal to a closed interval because each compact subset of $\R$ is contained in a closed finite interval.  Let $\phi^\alpha$ be a stationary optimal policy for a discount factor $\a\in [0,1)$ and $x^\a$ be a state such that $v_\alpha(x^\alpha)=v^{\phi^\alpha}_\a(x)=m_\alpha.$  Since $x^\alpha\in X_\a,$ then $x^\a\in [x^*_L,x^*_U].$ Consider a policy $\sigma$ such that, if the initial point $x<x_L,$ then $\sigma$ orders up to the level that the policy $\phi^\alpha$ would order at state $x^\alpha,$ and then $\sigma$ makes the same decisions as  $\phi^\alpha.$  Since a move from state $x_\tau$ to $x_{\tau+1}$ can be presented as two instant moves: from $x_\tau$ to $x^\a$ and from $x^\a$ to $x_{\tau+1}, $ as shown on Figure~\ref{verify_asmpB}, then
\begin{equation}\label{eq6.10} v_\alpha(x) \le v_\alpha^\sigma(x)\le K+\c (x^\alpha-x)+v_\alpha(x^\alpha)\le K+\c (x^*_U-x)+m_\alpha,\qquad\qquad x<x^*_L.
\end{equation}
For the initial inventory level $x\ge x^*_L,$  the policy $\sigma$ is defined in the following way.  It does not order as long as the inventory level is greater than or equal to $x_L^*.$  Then, as soon as the inventory level is less than $x^*_L,$ the policy $\sigma$ behaves in the same way as if it would behave if $x_\tau$ were the starting point, where    $\tau:=\inf\{\tau=0,1,\ldots: x_t<x^*_L\}$ is the first epoch when the inventory level is less than $x_L^*.$ Standard arguments from renewal theory imply that $\E [x_\tau]>-\infty$ and $C(x,\tau)<+\infty,$ where  $C(x,\tau)$ is the expected total undiscounted holding (or backordering) cost paid until the system reaches the level $x_\tau.$  Then
\begin{equation}\label{eq6.11} v_\alpha(x)\le v_\alpha^\sigma(x)\le
 C(x,\tau) +K +\c \E[x^*_U-x_\tau]+m_\alpha,\qquad\qquad x\ge x^*_L.
\end{equation}
Inequalities \eqref{eq6.10} and \eqref{eq6.11} imply that Assumption~{\bf B(ii)} holds. Though the above proof was applied in~\cite{FL1} to the classic periodic review system with backorders, it is generic and applicable to other systems.  For problems with lost sales the proof may be even simpler because it may be possible to define $\sigma$ so that $\tau$ is the first time when there is no inventory. Then  $x_\tau=0,$ and the expected cost of a lost sale will be added to the right hand side of \eqref{eq6.10}.  This expected cost is typically finite.

 Certain average cost
optimal policies can be approximated by discount optimal policies with vanishing
discount factor; see Feinberg et al.~\cite[Theorem 5]{FKZMDP}.  The following theorem and its corollary follow from such approximations.  In particular, the theorem and its corollary are useful for verifying that a limit point of optimal thresholds
for vanishing discount factors is an optimal threshold for average costs per unit time.

Recall that, for the function $u(x)$ defined in
\eqref{eq131}, for each $x\in\X$ there exist  sequences $\{\alpha_n\uparrow 1\}$
and $\{x^{(n)}\to x\},$ where $x^{(n)}\in\X,$ $n=1,2,\ldots,$ such that
$u(x)=\lim_{n\to\infty} u_{\alpha_n}(x^{(n)} ).$ Similarly, for a sequence
$\{\alpha_n\uparrow 1\}$ consider the function $\tilde u$ defined in
\eqref{eq1311}. Then for each $x\in\X$ there exist a sequence $\{x^{(n)}\to x\}$
of points in $\X$ and a subsequence $\{\alpha^*_n\}_{n=1,2,\ldots}$ of the
sequence $\{\alpha_n\}_{n=1,2,\ldots}$ such that $\tilde{u}(x)=\lim_{n\to\infty}
u_{\alpha^*_n}(x^{(n)} ).$
\begin{thm}\label{l:averoptt} (\cite{FL1}).
    Let Assumptions {\bf W*} and {\bf B} hold. For $x\in\X$ and $a^*\in\A,$ the
    following two statements hold:
    \begin{enumerate}[(i)]
     \item  for a sequence $\{( x^{(n)},\alpha_n)\}_{n=1,2,\ldots}$ with
          $0\le \alpha_n\uparrow 1,$ $ x^{(n)}\in\X,$   $ x^{(n)}\to x,$   and
          $u_{\alpha_n}(x^{(n)})\to u(x)$ as $n\to\infty,$ if there are a sequence
          of natural numbers $\{n_k\to\infty\}_{k=1,2,\ldots}$  and actions
          $\{a^{(n_k)}\in
          A_{\alpha_{n_k}}(x^{(n_k)})\}_{k=1,2,\ldots},$ such that $ a^{(n_k)}\to a^*$ as $k\to\infty,$
          then $a^* \in A^*_u(x),$ where the function $u$ is defined in \eqref{eq131};
          \label{state:converge1}
      \item let $\{\alpha_n\uparrow 1\}_{n=1,2\ldots}$ be   a sequence of discount factors, $\{\alpha^*_n\}_{n=1,2\ldots}$ be its subsequence, and   $\{x^{(n)}\to x\}_{n=1,2,\ldots}$ be a sequence of states from $\X$  such that $u_{\alpha^*_{n}}(x^{(n)})\to \tilde{u}(x)$ as
          $n\to\infty,$ where the function $\tilde u$ is defined in \eqref{eq1311} for
          the sequence $\{\alpha_{n}\}_{n=1,2,\ldots}.$
          If there are actions
          $a^{(n)}\in A_{\alpha^*_{n}}(x^{(n)})$ such that $ a^{(n)}\to a^*$ as
          $n\to\infty,$ then $a^*\in A^*_{\tilde u}(x).$ \label{state:converge2}
    \end{enumerate}
\end{thm}

\begin{cor}\label{l:averopt} (\cite{FL1}).
    Let Assumptions {\bf W*} and {\bf B} hold.  For $x\in\X$ and $a^*\in \A,$ the following two statements hold:
    \begin{enumerate}[(i)]
     \item if each sequence $\{(\alpha^*_n, x^{(n)})\}_{n=1,2,\ldots}$ with
          $0\le \alpha^*_n\uparrow 1,$ $x^{(n)}\in\X,$ and $ x^{(n)}\to x,$ $n=1,2,\ldots,$ contains a
          subsequence $(\alpha_{n_k}, x^{(n_k)}),$ such that there exist actions
          $a^{(n_k)}\in A_{\alpha_{n_k}}(x^{(n_k)})$ satisfying $ a^{(n_k)}\to
          a^*$ as $k\to\infty,$ then $a\in A^*_u(x)$ with the function $u$ defined
          in \eqref{eq131}; \label{state:averopt1}
      \item  if there is a sequence $\{\alpha_n\uparrow 1\}_{n=1,2,\ldots},$ such that for every sequence of states $\{x_n\to x\}$ from $\X$ there are actions $a^{n}\in A_{\alpha_n}(x^{(x)}),$ $n=1,2,\ldots,$   satisfying $a_n\to a^*$ as $n\to\infty,$
          then $a^*\in A^*_{\tilde u}(x),$  where the function $\tilde u$ is defined in \eqref{eq1311} for
          the sequence $\{\alpha_{n}\}_{n=1,2,\ldots}.$ \label{state:averopt2}
    \end{enumerate}
\end{cor}

The following theorem is useful for proving asymptotic properties of optimal actions for discounted problems when  the discount factor tends to 1.

\begin{thm}\label{l:averopttt} (\cite{FL1}).
    Let Assumptions {\bf W*} and {\bf B} hold. For $x\in\X$  the
    following two statements hold:
    \begin{enumerate}[(i)]
    \item there exists a compact set $D^*(x)\subseteq \A$ such that  $ A_\alpha(x)\subseteq D^*(x)$ for all $\alpha\in[0,1);$
    \item if $\{\alpha_n\}_{n=1,2\ldots}$ is a sequence of discount factors  $ \alpha_n\in [0,1),$  then every sequence of infinite-horizon $\alpha_n$-optimal actions $\{a^{(n)}\in A_{\alpha_n}(x)\}_{n=1,2,\ldots}$  is bounded and therefore  has a limit point  $a^*\in \A.$
        \end{enumerate}
    \end{thm}

\section{Partially Observable Markov Decision Processes}\label{s6}
POMDPs model the situations,
when the current state of the system may be unknown, and the decision maker
uses indirect observations for decision making.  A POMDP is defined by the same objects as an MDP, but in addition to the state space
$\X$ and action space $\A.$   The states and observations are linked by the transition probability $Q(dy_{t+1}|a_t,x_{t+1}),$ from $\A\times\X$ to $\Y,$ $t=0,1,\ldots.$  Thus, a POMDP is defined as the tuple $\{\X,\Y,\A, P,Q,c\},$ where
the Borel state and action spaces $\X$ and $\A,$  the transition probability $P, $ and the cost function $c$ are the same objects as in an MDP.  In addition,  $\Y$ is the observation space, which is also assumed to be a Borel subset of a Polish space, and $Q$ is the observation probability, which is a regular transition probability from $\A\times\X$ to $\Y.$  Sometimes we say a transition kernel or a stochastic kernel instead of transition probability.  Though the initial state of the system may be unknown, the decision maker knows the probability distribution of the initial state $p(dx_0),$ and there is an observation probability for the first observation $Q_0(dy_0|x_0).$

  In various applications it is possible that there are continuous states and discrete observations, discrete states and continuous observations, and  both spaces can be discrete or continuous.  So, we consider a general situation by assuming that $\X$ and $\Y$ are Borel subsets of Polish spaces.

The following subsection describes a classic transformation of a POMDP to a Completely Observable MDP (COMDP), whose states are posterior probability distributions of states in the POMDP.  This transformation was introduced by Aoki~\cite{Ao}, \AA str\"om~\cite{As},  Dynkin~\cite{Dyn}, and Shiryaev~\cite{sh}.  These ideas were advanced in the book by Striebel~\cite{Str} and in the references provided in the following subsection. The main results of this section describe optimality conditions for POMDPs and COMDPs introduced in Feinberg et al.~\cite{FKZg}.

\begin{figure}[h]
  \centering
  \caption{POMDP Diagram.}
  \label{POMDP}
  \includegraphics[scale=0.12]{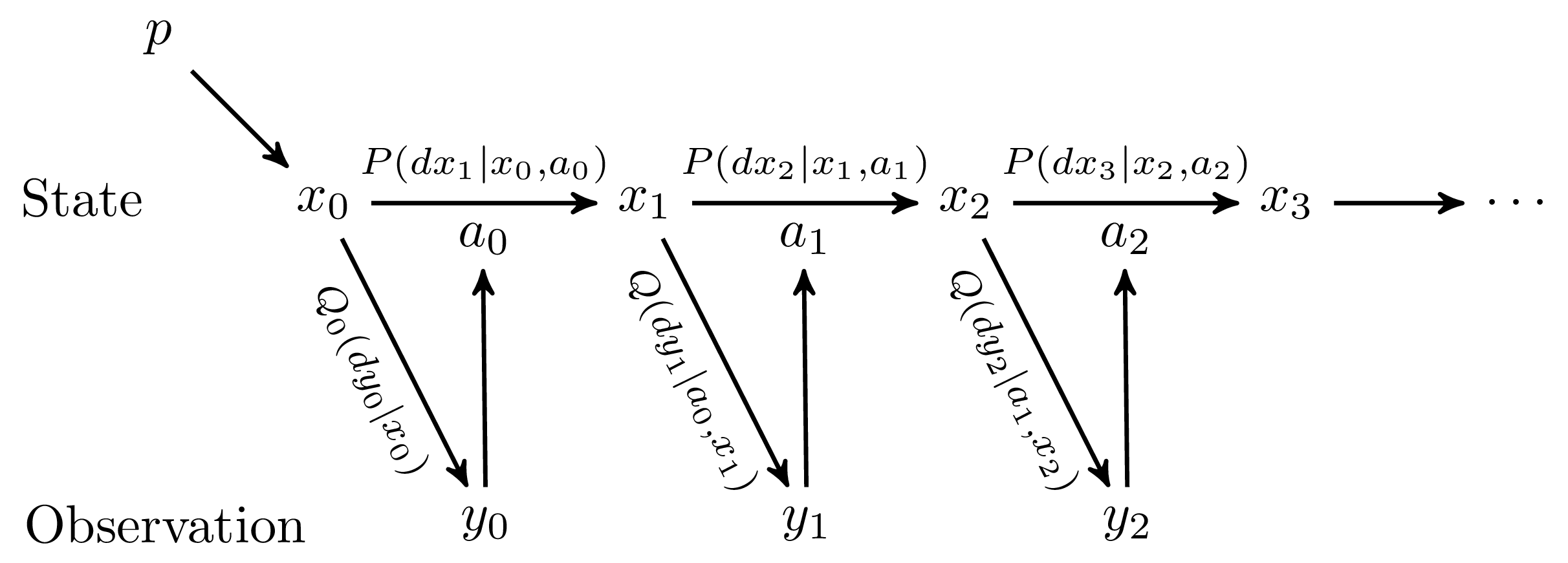}
\end{figure}

The POMDP evolves as follows.
At time $t=0$, the initial unobservable state $x_0$ has a given
prior distribution $p.$ The initial observation $y_0$ is generated
according to the initial observation kernel $Q_0(\,\cdot\,|x_0).$ At
each time epoch $\n=0,1,\ldots,$ if the state of the system is
$x_\n\in\X$ and the decision-maker chooses an action $a_\n\in \A$,
then the cost $c(x_\n,a_\n)$ is incurred;
 the system moves to state $x_{\n+1}$ according to
the transition law $P(\,\cdot\,|x_\n,a_\n).$ The observation
$y_{\n+1}\in\Y$ is generated by the observation kernels
$Q(\,\cdot\,|a_\n,x_{\n+1})$, $\n=0,1,\ldots,$ and
$Q_0(\,\cdot\,|x_0);$ see Figure~\ref{POMDP}. For the state space $\X$, denote by
$\P(\X)$ the set of probability measures on $(\X,{\cal B}(\X)).$  We always consider a metric on $\P(\X)$
consistent with the topology of weak convergence.

Define the \textit{observable histories}: $h_0:=(p,y_0)\in H_0$
and
$h_\n:=(p,y_0,a_0,\ldots,y_{n-1}, a_{\n-1}, y_\n)\in H_\n$ for all
$n=1,2,\dots,$
where $H_0:=\P(\X)\times \Y$ and $H_\n:=H_{\n-1}\times \A\times \Y$
if $\n=1,2,\dots\ .$ Then a \textit{policy} for the POMDP is defined as
a sequence $\pi=\{\pi_\n\}$ such that, for each $n=0,1,\dots,$
$\pi_\n$ is a transition kernel on $\A$ given $H_\n$.  Moreover, $\pi$
is called \textit{nonrandomized}, if each probability measure
$\pi_\n(\cdot|h_\n)$ is concentrated at one point.  The \textit{set of all policies} is denoted by $\Pi$.
The Ionescu Tulcea theorem (Bertsekas and Shreve \cite[pp.
140-141]{besh96} or Hern\'andez-Lerma and Lassere
\cite[p.178]{HL:96}) implies that, given a policy $\pi\in \Pi,$  an
initial distribution $p\in \P(\X)$ and a sequence of  transition probabilities
  $Q_0,\pi_0,P,Q,\pi_1,P,Q,\pi_2,\ldots$ determine a unique probability measure
$P_{p}^\pi$ on the set of all trajectories
$\mathbb{H}_{\infty}=(\X\times\Y\times \mathbb{A})^{\infty}$
endowed with the $\sigma$-field, which is the product of  Borel
$\sigma$-fields on $\X$, $\Y$, and $\mathbb{A}$ respectively.
The expectation with respect to this probability measure is
denoted by $\E_{p}^\pi$.

Let us specify a performance criterion. For a finite horizon
$N=0,1,\ldots,$ and for a policy $\pi\in\Pi$, let  the
\textit{expected total discounted costs} be
\begin{equation}\label{eq1}
v_{N,\alpha}^{\pi}(p):=\mathbb{E}_p^{\pi}\sum\limits_{\n=0}^{N-1}\alpha^\n c(x_\n,a_\n),\
p\in \P(\X),
\end{equation}
where $\alpha\ge 0$ is the discount factor, and
$v_{0,\alpha}^{\pi}(p)=0.$ When $N=\infty$, we always assume $\alpha\in [0,1).$  We always assume that the
function $c$ is bounded below.

For any function $g^{\pi}(p)$, including
$g^{\pi}(p)=v_{N,\alpha}^{\pi}(p)$ and
$g^{\pi}(p)=v_{\alpha}^{\pi}(p)$ define the \textit{optimal cost}
\begin{equation*} g(p):=\inf\limits_{\pi\in \Pi}g^{\pi}(p), \qquad
\ p\in\P(\X),
\end{equation*} where $\Pi$ is \textit{the set of all policies}.
A policy $\pi$ is called \textit{optimal} for the respective
criterion, if $g^{\pi}(p)=g(p)$ for all $p\in \P(\X).$ For
$g^\pi=v_{N,\alpha}^\pi$, the optimal policy is called
\emph{$N$-horizon discount-optimal}; for $g^\pi=v_{\alpha}^\pi$, it
is called \emph{discount-optimal}.

\subsection{Reduction of POMDPs to MDPs}

In this section, we 
%
formulate the well-known reduction of a POMDP to the corresponding
COMDP (\cite{besh96, DY, HLerma1, Rh, Yu}).  This reduction constructs an MDP whose states are
probability distributions on the original state space.  These distributions are posteriori distributions of states
after the observations become known.  In addition to posterior probabilities, they are also called belief probabilities and belief states in the literature.  The reduction establishes the correspondence between certain classes of policies in MDPs and POMDPs and their performances.  If an optimal policy is found for the COMDP, it defines in a natural way an optimal policy for the original POMDP.  The reduction holds for measurable transition probabilities, observation probabilities, and one-step costs.  Except for problems with discrete transition probabilities or with transition probabilities having densities (see \cite{BR, Bens2}), almost nothing had been known until recently on the existence of optimal policies for POMDPs and how to find them.

To simplify notations,
we sometimes drop the time parameter.  Given a posterior
distribution $z$ of the state $x$ at time epoch  $\n=0,1,\ldots$
and given an action $a$ selected at epoch $\n$, denote by
$R(B\times C|z,a) $ the joint probability that the  state at time
$(\n+1)$  belongs to the set $B\in {\mathcal B}(\X)$ and the
observation at time $(\n+1)$ belongs to the set $C\in {\mathcal
B}(\Y)$,
\begin{equation}\label{3.3}
R(B\times
C|z,a):=\int\limits_{\X}\int\limits_{B}Q(C|a,x')P(dx'|x,a)z(dx),
\end{equation}
where $R$ is a transition kernel on $\X\times\Y$ given
${\P}(\X)\times \A$;  see Bertsekas and Shreve \cite{besh96},
Dynkin and Yushkevich \cite{DY},  Hern\'{a}ndez-Lerma
\cite{HLerma1}, or Yushkevich \cite{Yu} for details.  Therefore,
 the probability $R'( C|z,a) $ that the observation $y$ at time $\n$
belongs to the set $C\in \mathcal{B}(\Y)$ is
\begin{equation}\label{3.5}
R'(C|z,a)=\int\limits_{\X}\int\limits_{\X}Q(C|a,x')P(dx'|x,a)z(dx),
\end{equation}where $R'$ is a transition kernel on $\Y$ given
${\P}(\X)\times \A.$ By Bertsekas and Shreve~\cite[Proposition
7.27]{besh96}, there exists a transition kernel $H$ on $\X$ given
${\P}(\X)\times \A\times\Y$ such that
\begin{equation}\label{3.4}
R(B\times C|z,a)=\int\limits_{C}H(B|z,a,y)R'(dy|z,a),
\end{equation}

The transition kernel $H(\,\cdot\,|z,a,y)$ defines a measurable mapping $H:\,\P(\X)\times \A\times \Y \to\P(\X)$, where
$H(z,a,y)[\,\cdot\,]=H(\,\cdot\,|z,a,y).$ For each pair $(z,a)\in \P(\X)\times\A$, the mapping $H(z,a,\cdot):\Y\to\P(\Y)$ is defined
$R'(\,\cdot\,|z,a)$-a.s. uniquely in $y$; Dynkin and Yushkevich \cite[p.~309]{DY}. It is known that for a posterior distribution $z_\n\in \P(\X)$, action
$a_\n\in A(x)$, and an observation $y_{\n+1}\in\Y,$ the posterior distribution $z_{\n+1}\in\P(\X)$ is
\begin{equation}\label{3.1}
z_{\n+1}=H(z_\n,a_\n,y_{\n+1}).
\end{equation}
However, the observation $y_{\\n+1}$ is not available in the COMDP
model, and therefore $y_{\n+1}$ is a random variable with the
distribution  $R'(\,\cdot\,|z_\n,a_\n)$, and (\ref{3.1}) is a
stochastic equation that maps $(z_\n,a_\n)\in \P(\X)\times\A$ to
$\P(\X).$ The stochastic kernel that defines the distribution of
$z_{\n+1}$ on $\P(\X)$ given $\P(\X)\times\X$ is defined uniquely as
\begin{equation}\label{3.7}
q(D|z,a):=\int\limits_{\Y}\h_D[H(z,a,y)]R'(dy|z,a),
\end{equation}
where for $D\in \mathcal{B}(\P(\X))$
\[
\h_D[u]=\left\{
\begin{array}{ll}
1,&u\in D,\\
0,&u\notin D;
\end{array}
\right.
\]
Hern\'andez-Lerma~\cite[p. 87]{HLerma1}. The measurable particular
choice of stochastic kernel $H$ from (\ref{3.4}) does not affect on
the definition of $q$ from (\ref{3.7}), since for each pair
$(z,a)\in \P(\X)\times\A$, the mapping $H(z,a,\cdot):\Y\to\P(\Y)$ is
defined $R'(\,\cdot\,|z,a)$-a.s. uniquely in $y$; Dynkin and
Yushkevich \cite[p.~309]{DY}.

The COMDP is defined as an MDP with parameters
($\P(\X)$,$\A$,$q$,$\c$), where
\begin{itemize}
\item[(i)] $\P(\X)$ is the state space; \item[(ii)] $\A$ is the
action set available at all states $z\in\P(\X)$; \item[(iii)] the
 one-step cost function $\c:\P(\X)\times\A\to\R$, defined as
\begin{equation}\label{eq:c}
\c(z,a):=\int\limits_{\X}c(x,a)z(dx), \quad z\in\P(\X),\, a\in\A;
\end{equation}
 \item[(iv)] the transition probabilities $q$ on $\P(\X)$
given $\P(\X)\times \A$ defined in (\ref{3.7}).
\end{itemize}

If a stationary optimal policy for the COMDP exists and is found, it
allows the decision maker to formulate an optimal policy for the
POMDP.  The details on how to do this can be found in Bertsekas and Shreve~\cite{besh96} or
Dynkin and Yushkevich~\cite{DY},  or Hern\'{a}ndez-Lerma~\cite{HLerma1}.  Therefore,
a POMDP can be reduced to a COMDP.
This reduction holds for measurable transition kernels $P$, $Q$,
$Q_0$. The measurability of these kernels and the cost function $c$
lead to the measurability of transition probabilities for the
corresponding COMDP.

As follows from Theorem~\ref{prop:dcoe}, if the COMDP satisfies Assumption~{\bf W*}, then optimal policies exist, they satisfy the optimality equation, and can be found by value iterations.   This is formulated in Theorem~\ref{teor4.3a} below. The validity of Assumption~{\bf W*} for the COMDP is equivalent to the correctness of the following two Hypotheses:

\textbf{Hypothesis (i).} The transition probability $q$ from $\P(\X)\times\A$ to $\P(\X))$ is weakly continuous.

\textbf{Hypothesis (ii).} The cost function $\c:\P(\X)\times\A\to {\bar\R}$  is bounded below and $\K$-inf-compact on $\P(\X)\times\A.$

Following theorem states the correctness of  Hypothesis (ii).  The question, whether Hypothesis (i) holds, a more difficult, and the the following subsection is devoted to answering it.

\begin{thm}\label{th:wstar} (\cite{FKZg}).
If the function $c:\X\times\A\to \R$ is a bounded  below,
$\K$-inf-compact (inf-compact) function on $\X\times\A$, then the cost function
$\c:\P(\X)\times\A\to\R$ defined for the COMDP in (\ref{eq:c}) is
bounded  below by the same constant and $\K$-inf-compact (inf-compact)  on $\P(\X)\times\A$.
\end{thm}

In addition to weak convergence, two types of convergence are mentioned in the next subsection: setwise convergence and  convergence in total variation. Here we recall their definitions.

 Let $(P_n)_{n=1,2,\ldots}$  be a sequence of probability measures on a measurable space $(\S,\cal F)$.  This sequence converges setwise to a probability measure $P_0$ on  $(\S,\cal F)$ if $\lim_{n\to\infty} P_n(A)=P_0(A)$ for each $A\in \cal F.$  This sequence converges in total variation if  $\lim_{n\to\infty} ||P_n-P_0||=0, $ where  $||P_n(A)-P_0(A)||=2\sup\{P_n(A)-P_0(A): A\in\cal F\}.$  Convergence in total variation implies setwise convergence.  If $\S$ is a metric space and $\cal F$ is its Borel $\sigma$-field, then setwise convergence implies weak convergence. Recall that  $P^*$ is a regular transition probability from a metric space $\S_1$ to a metric space $\S_2,$ if $P^*(\cdot|s_1)$ is a probability measure on $\S_1$ for each $s\in \S_2$ and $P^*(A|\cdot)$ is a Borel function on $\S_1$ for each Borel subset $A$ of $\S_2.$ A transition probability is weakly (setwise, in total variation) continuous, if, for every sequence $(s^n)_{n=1,2\ldots}$ on $S_1$ converging to $s^0\in\S_1,$
the sequence $(P^*(\cdot|s^n))_{n=1,2,\ldots}$ converges weakly (setwise, in total variation)
to $P^*(\cdot|s^0).$  There are two mathematical tools that are useful for the analysis of convergence of probability measures and for the analysis of MDPs and POMDPs: Fatou's lemma for variable probabilities (see Feinberg et al.~\cite{FKZTV} and references therein) and uniform Fatou's lemma introduced in Feinberg et al.~\cite{FKZUFL}.
\subsection{Optimality Conditions for Discounted POMDPs}

For the COMDP, Assumption~{\bf W*} can be rewritten in the
following form:

(i)  $\c$ is $\K$-inf-compact on $\P(\X)\times\A$;


(ii) the transition probability  $q(\cdot|z,a)$ is weakly continuous in $(z,a)\in \P(\X)\times\A$.


Theorem~\ref{prop:dcoe}  has the following form for the COMDP
$(\P(\X),\A,q,\c)$:

\begin{thm}{\rm (cf. Feinberg et al. \cite[Theorem~2]{FKZMDP}).} \label{teor4.3a}
Let the COMDP $(\P(\X),\A,q,\c)$ satisfy Assumption~{\bf
W*}. 
Then:

{(i}) the functions $v_{\n,\alpha}$, $\n=0,1,\ldots$, and $v_\alpha$
are lower semi-continuous on $\P(\X)$, and
$v_{\n,\alpha}(z)\to 
 v_\alpha (z)$ as $\n \to \infty$ for all
$z\in \P(\X);$

{(ii)} for any $z\in \P(\mathbb{X})$, and $\n=0,1,...,$
\begin{equation}\label{eq433}\begin{aligned}
&v_{\n+1,\alpha}(z)=\min\limits_{a\in \A}\{\c(z,a)+\alpha
\int_{\P(\X)} v_{\n,\alpha}(z')q(dz'|z,a)\}\\ & =
\min\limits_{a\in
\A}\{\int\limits_{\X}c(x,a)z(dx)+\int\limits_{\X}\int\limits_{\X}\int\limits_{\Y}
v_{\n,\alpha}(H(z,a,y))\\&\times\alpha  Q(dy|a,x')P(dx'|x,a)z(dx)
\},\end{aligned}
\end{equation}
where $v_{0,\alpha}(z)=0$ for all $z\in \P(\X)$, and the nonempty sets
\[\begin{aligned}
&A_{\n,\alpha}(z):=\{ a\in \A:\,v_{\n+1,\alpha}(z)\\ & =c(z,a)
 +\alpha \int_{\P(\X)} v_{\n,\alpha}(z')q(dz'|z,a)
\},
\end{aligned}\] where $z\in \P(\X)\ $, satisfy the following properties: (a) the graph
${\rm Gr}(A_{\n,\alpha})=\{(z,a):\, z\in\P(\X), a\in
A_{\n,\alpha}(z)\}$, $\n=0,1,\ldots,$ is a Borel subset of
$\P(\X)\times \mathbb{A}$, and (b) if $v_{\n+1,\alpha}(z)=\infty$,
then $A_{\n,\alpha}(z)=\A$ and, if $v_{\n+1,\alpha}(z)<\infty$, then
$A_{\n,\alpha}(z)$ is compact;

{(iii)} for any $N=1,2,\ldots$, there exists a Markov optimal
$N$-horizon  policy $(\phi_0,\ldots,\phi_{N-1})$ for the COMDP, and if for an
$N$-horizon Markov policy $(\phi_0,\ldots,\phi_{N-1})$  the
inclusions $\phi_{N-1-\n}(z)\in A_{\n,\alpha}(z)$, $z\in\P(\X),$
$\n=0,\ldots,N-1,$ hold, then this policy is $N$-horizon optimal;

{(iv)} for $\alpha\in [0,1)$
\[\begin{aligned}
 & v_{\alpha}(z)=\min\limits_{a\in
\A}\{\c(z,a)+\alpha\int_{\P(\X)}
v_{\alpha}(z')q(dz'|z,a)\}\\ & =  
%
\min\limits_{a\in \A}\{\int\limits_{\X}c(x,a)z(dx) +\alpha
\int\limits_{\X}\int\limits_{\X}\int\limits_{\Y}
v_{\alpha}(H(z,a,y))\\ & \times Q(dy|a,x')P(dx'|x,a)z(dx) \},\quad
z\in \P(\X), \end{aligned}\] and the nonempty sets
\[\begin{aligned}
A_{\alpha}(z):= & \{a\in \A:\,v_{\alpha}(z)=\c(z,a)\\
+ & \alpha\int_{\P(\X)} v_{\alpha}(z')q(dz'|z,a) \},\quad z\in
\P(\X),
\end{aligned}\] satisfy the following properties: (a) the graph
${\rm Gr}(A_{\alpha})=\{(z,a):\, z\in\P(\X), a\in \A_\alpha(z)\}$
is a Borel subset of $\P(\X)\times \mathbb{A}$, and (b) if
$v_{\alpha}(z)=\infty$, then $A_{\alpha}(z)=\A$ and, if
$v_{\alpha}(z)<\infty$, then $A_{\alpha}(z)$ is compact.

{(v)} for an infinite horizon there exists a stationary discount-optimal policy $\phi_\alpha$ for the COMDP, and a stationary policy $\phi$ is optimal if and only if
$\phi_\alpha(z)\in A_\alpha(z)$ for all $z\in \P(\X).$

{(vi)}  if the function $c$ is inf-compact,  the functions
$v_{\n,\alpha}$, $\n=1,2,\ldots$, and $v_\alpha$ are inf-compact on
$\P(\X)$.
\end{thm}


Hern\'andez-Lerma~\cite[Section 4.4]{HLerma1} provided the
following conditions for the existence of optimal policies for the
COMDP: (a) $\A$ is compact, (b) the cost function $c$ is bounded
and continuous, (c) the transition probability $P(\cdot|x,a)$ and
the observation kernel $Q(\cdot|a,x)$ are weakly continuous
transition kernels; (d) there exists a weakly continuous
$H:\P(\X)\times\A\times\Y\to\P(\X)$ satisfying (\ref{3.4}).
Consider the following relaxed version of assumption (d).

\noindent\textbf{Assumption~{\bf H}.} (\cite{FKZg}). There exists a transition kernel $H$ on $\X$ given $\P(\X)\times\A\times\Y$ satisfying (\ref{3.4}) such that: if a
sequence $\{z^n\}\subseteq\P(\X)$ converges weakly to $z\in\P(\X)$, and $\{a^n\}\subseteq\A$ converges to $a\in\A$, $n\to\infty$, then there exists a
subsequence $\{(z^{n_k},a^{n_k})\}_{k\ge 1}\subseteq \{(z^{n},a^{n})\}_{n\ge 1}$ such that
\[
H(z^{n_k},a^{n_k},y)\mbox{ converges weakly to }H(z,a,y),\
n\to\infty,
\]
and this convergence takes place  $R'(\,\cdot\,|z,a)$ almost
surely in $y\in\Y$.


%


The following theorem provides two sufficient conditions for weak continuity of $q.$  Statement (ii) can be found in Hernandez-Lerma~\cite[p. 90]{HLerma1}.
\begin{thm}\label{mainN} (\cite{FKZg}).  If the transition probability $P(dx'|x,a)$ is weakly continuous, then  each of the following two conditions implies weak continuity of the transition probability $q$ from $\P(\X)\times\A$ to $\P(\X):$
 \begin{itemize}

\item[(i)] the transition probability $R'(dy|z,a)$  from
$\P(\X)\times\A$ to $\Y$ is setwise continuous, and Assumption~{\bf H}
holds,

\item[(ii)]
 the transition probability $Q(dy|a,x)$  from
$\A\times\X$  to $\Y$ is weakly continuous, and there exists a weakly
continuous $H:\P(\X)\times\A\times\Y\to\P(\X)$ satisfying
(\ref{3.4}).\end{itemize}
\end{thm}

Weak continuity of the transition probability $P$ and continuity of the  transition probability $Q$  in  total variation imply that Assumption~{\bf H}
holds, and this leads to the following theorem.

\begin{thm}\label{teor:H} (\cite{FKZg}).
Let the transition probability $P(dx'|x,a)$   from $\X\times\A$ to $\X$
be weakly continuous and let the transition probability $Q(dy|a,x)$  from $\A\times\X$ to
$\Y$   be continuous in  total variation.
Then the transition probability $R'(dy|z,a)$  from
$\P(\X)\times\A$ to  $\Y$  is setwise continuous, Assumption~{\bf H}
holds, and  the transition probability $q$  from
$\P(\X)\times\A$  to  $\P(\X)$ is weakly continuous.
\end{thm}

The following theorem, which follows from Theorems~\ref{th:wstar}--\ref{mainN},   relaxes  assumptions (a), (b), and (d) in
Hern\'andez-Lerma~\cite[Section 4.4]{HLerma1}.

\begin{thm}\label{main} (\cite{FKZg}).
Under the following  conditions: \begin{itemize}
\item[(a)] the cost function $c$ is  $\K$-inf-compact; 
\item[(b)] either
\begin{enumerate}[{(i)}]
\item the transition probability $R'(dy|z,a)$  from
$\P(\X)\times\A$ to $\Y$  is setwise continuous and Assumption~{\bf H}
holds,

or \item the transition probability $Q(dy|a,x)$ from
$\A\times\X$  to $\Y$ is weakly continuous and there exists a weakly
continuous $H:\P(\X)\times\A\times\Y\to\P(\X)$ satisfying
(\ref{3.4});\end{enumerate}\end{itemize}
the COMDP $(\P(\X),\A,q,\c)$ satisfies Assumption~{\bf W*} and therefore statements (i)--(vi) of
Theorem~\ref{teor4.3a} hold.
\end{thm}
%

Theorems~\ref{teor:H} and \ref{main}  imply the following result.

\begin{thm}\label{cteor:H} (\cite{FKZg}).
Let Assumption~{\bf W*} hold and
let the transition probability $Q(dy|a,x)$ from $\A\times\X$ to $\Y$
be continuous in total variation. Then  statements (i)--(vi)
of Theorem~\ref{teor4.3a} hold.
\end{thm}

Theorem~\ref{main} assumes either the weak continuity of $H$ or Assumption~{\bf H} together with the setwise continuity of $R'$. For some
applications, including the inventory control applications described in Section~\ref{s7}, the filtering kernel $H$ satisfies Assumption~{\bf H} for some
observations and it is weakly continuous for other observations. The following theorem is applicable to such situations.

\begin{thm}\label{cor:main} (\cite{FKZg}).
Let the observation space $\Y$ be partitioned into two disjoint
subsets $\Y_1$ and $\Y_2$ such that $\Y_1$ is open in $\Y$. Suppose the
following assumptions hold:

(a) the transition probabilities $P$ to $\X$ from $\X\times\A$  to $\X$ and
$Q$  from $\A\times\X$ to $\Y$ are weakly continuous;

(b) the measure $R'(\,\cdot\,|z,a)$ on $(\Y_2,\B(\Y_2))$ is setwise
continuous in $(z,a)\in \P(\X)\times \A,$ that is, for every
sequence $\{(z^n,a^n)\}_{n=1,2,\ldots}$  in $\P(\X)\times\A$
converging to $(z,a)\in\P(\X)\times\A$  and for every
$C\in\B(\Y_2),$ we have
$R'(C|z^n,a^n)\to R'(C|z,a);$ 

 (c)  there exists a transition probability
$H$  from $\P(\X)\times\A\times\Y$ to $\X$ satisfying (\ref{3.4})
such that: \begin{enumerate}[{(i)}]
\item the transition probability $H$  from
$\P(\X)\times\A\times\Y_1$  to $\X$  is weakly continuous;
%
 \item  Assumption~{\bf H} holds on $\Y_2,$ that is, if a
sequence $\{z^{(n)}\}_{n=1,2,\ldots}\subseteq\P(\X)$ converges
weakly to $z\in\P(\X)$ and a sequence
$\{a^{(n)}\}_{n=1,2,\ldots}\subseteq\A$ converges to $a\in\A$,
then there exists a subsequence $\{(z^{(n_k)},a^{(n_k)})\}_{k=1,2,\ldots}\subseteq
\{(z^{(n)},a^{(n)})\}_{n=1,2,\ldots}$ and a measurable subset $C$ of $\Y_2$ such that
$R'(\Y_2\setminus C|z,a)=0$ and  
$H(z^{(n_k)},a^{(n_k)},y)$  converges weakly to
$H(z,a,y)$ for all $y\in C$; 
\end{enumerate}

\noindent
Then the transition probability $q$ from $\P(\X)\times\A$ to $\P(\X)$
is weakly continuous.  If, in addition to the above conditions,
 the cost function $c$ is  $\K$-inf-compact, then the
COMDP $(\P(\X),\A,q,\c)$ satisfies Assumption~{\bf W*} and therefore statements (i)--(vi) of
Theorem~\ref{teor4.3a} hold.

\end{thm}

The following corollary follows from Theorem~\ref{cor:main}.

\begin{cor}\label{cor:main1} (\cite{FKZg}).
Let the observation space $\Y$ be partitioned into two
disjoint subsets $\Y_1$ and $\Y_2$ such that $\Y_1$ is  open in $\Y$ and $\Y_2$ is countable. Suppose the
following assumptions hold:

(a) the transition probabilities $P$ from $\X\times\A$ to $\X$  and
$Q$  from $\A\times\X$ to $\Y$ are weakly continuous;

(b)
$Q(y|a,x)$ is a continuous function on $\A\times\X$ for each $y\in
\Y_2;$

(c) there exists a stochastic kernel $H$ on $\X$ given
$\P(\X)\times\A\times\Y$ satisfying (\ref{3.4}) such that the
stochastic kernel $H$ on $\X$ given $\P(\X)\times\A\times\Y_1$ is
weakly continuous.\\
Then assumption (b) and (ii) from Theorem~\ref{cor:main} hold  and
the transition probability $q$ from $\P(\X)\times\A$  to $\P(\X)$ is
weakly continuous.  If, in addition to the above conditions,
the cost function $c$ is  $\K$-inf-compact, then the COMDP
$(\P(\X),\A,q,\c)$ satisfies Assumption~{\bf W*}
and therefore statements (i)--(vi) of Theorem~\ref{teor4.3a} hold.
\end{cor}

In conclusion of this section, we would like to mention another model of a controlled Markov process
with partial observations, in which the observation kernel $Q$ is not
defined explicitly, and a state of the system consists of two parts:
one part of the state is observable and another one is not; see e.g.,
Rhenius~\cite{Rh}, Yushkevich\cite{Yu}, B\"auerle and
Rieder\cite[Chapter 5]{BR}.  In Feinberg et al.~\cite{FKZST, FKZg}
such models were called  Markov Decision Models with Incomplete
Information, and the most general known sufficient conditions for the
existence of optimal policies for such models with the expected total
costs are provided in Feinberg et al~\cite[Theorem 6.2]{FKZST}.

\section{Inventory Control with Incomplete Information}\label{s7}
Bensoussan et al. \cite{Bens3}--\cite{Bens5} studied several inventory  control  problems for periodic review systems,
when the Inventory Manager (IM) may not have complete information
about inventory levels. In Bensoussan et
al.~\cite{Bens3}, \cite{Bens5}, a  problem with backorders is
considered.  In the model considered in \cite{Bens3}, the IM does
not know the inventory level, if it is nonnegative, and the IM
knows the inventory level, if it is negative. In the model
considered in \cite{Bens5},  the IM only knows  whether the
inventory level is negative or nonnegative.  In \cite{Bens4} a
problem with lost sales is studied where the IM only knows whether
a lost sale happened or not. The underlying mathematical analysis
is summarized in \cite{BensOx}, where additional references can be
found.  The analysis includes transformations of density functions
of demand distributions.

This section describes periodic review systems with backorders and lost sales, when some inventory levels are observable and some are not.  The goal
is to minimize the expected total costs.  Demand distributions  may not have densities.  This model is introduced in~Feinberg et al. \cite[Section 8.2]{FKZg}.

In the case of full observations, we model the problem as an MDP
with the state space $\X=\mathbb{R}$ (the current inventory level),
action space $\A=\mathbb{R}$ (the ordered amount of inventory),
and action sets $\A(x)=\A$ available at states $x\in \X$. If in a
state $x$ the amount of inventory $a$ is ordered, then the
holding/backordering cost $h(x)$,  ordering cost $C(a),$ and lost
sale cost $G(x,a)$ are incurred, where it is assumed that $h,$
$C,$ and $G$ are
nonnegative lower semi-continuous functions 
with values in ${\R}$ and $C(a)\to +\infty$ as $|a|\to \infty.$ Observe that the one-step cost function $c(x,a)=h(x)+C(a)+G(x,a)$ is $\K$-inf-compact on
$\X \times \A$.  For problems with back orders (no lost sales), usually $G(x,a)=0$ for all $x$ and $a.$

Let $D_t, t = 1,2,\ldots,$ be i.i.d. random variables with the distribution function $F_D$, where $D_t$ is the demand at epoch $t.$  The
dynamics of the system are defined by $x_{t+1}=F(x_t,a_t,D_{t+1}),$ where $x_t$ is the current inventory level  and $a_t$ is the ordered (or scrapped) inventory
at epoch
$t=0,1,\ldots\ .$ 
For problems with backorders $F(x_t,a_t,D_{t+1})=x_t+a_t-D_{t+1}$ and for
problems with lost sales $F(x_t,a_t,D_{t+1})=(x_t+a_t-D_{t+1})^+$.  In
both cases, $F$ is a continuous function defined on
$\mathbb{R}^3$.  To simplify and unify the presentation, we do not
assume
 $\X=[0,\infty)$
for models with lost sales. 
However,
 for problems with lost sales it is assumed that the initial state distribution $p$ is concentrated on $[0,\infty)$, and this implies that states $x<0$ will never
be visited. We assume that the distribution function $F_D$ is atomless (an equivalent assumption is that the function $F_D$ is continuous).
The state transition law $P$ on $\X$ given $\X\times\A$ is
\begin{equation}
\label{eq:P} P(B|x,a)=\int_{\mathbb{R}}\h\{F(x,a,s)\in B\}dF_D(s),
\end{equation}
where $B\in \B(\X),$ $x\in\X,$ and  $a\in\A.$
Since we do not assume that demands are nonnegative, this model also covers cash balancing problems and problems with returns; see Feinberg and
Lewis~\cite{FL} and the references therein.  In a particular case, when $C(a)=+\infty$ for $a<0$,  orders with negative sizes are infeasible, and, if an order
is placed, the ordered amount of inventory should be positive.

\begin{figure}[h]
  \centering
  \caption{Example with known current inventory level.}
  \label{inventory_knownLevel}
  \includegraphics[scale=0.3]{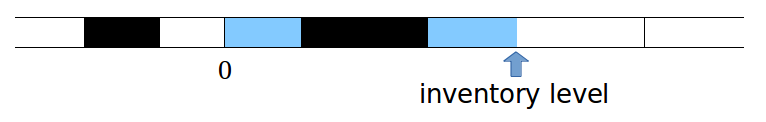}
\end{figure}

\begin{figure}[h]
  \centering
  \caption{Example with unknown current inventory level.}
  \label{inventory_unknownLevel}
  \includegraphics[scale=0.3]{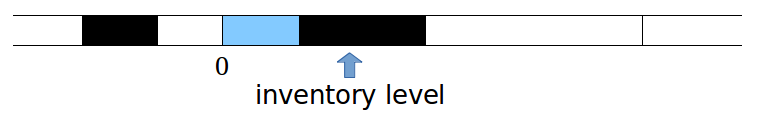}
\end{figure}

\begin{figure}[h]
  \centering
  \caption{Example with known inventory level and unknown backorder level.}
  \label{inventory_unknownBackorders}
  \includegraphics[scale=0.3]{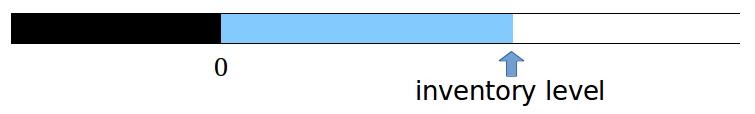}
\end{figure}

\begin{figure}[h]
  \centering
  \caption{Example with inventory level 0 and current inventory level inside a (nontransparent) container.}
  \label{inventory_unknownZeroLv}
  \includegraphics[scale=0.3]{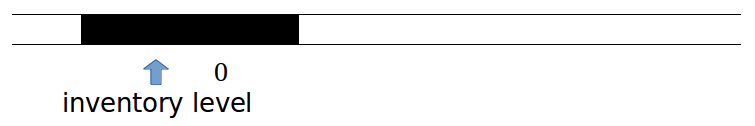}
\end{figure}

As mentioned above, some states (inventory levels) $x\in\X=\mathbb{R}$ are observable and some are not.  Let the inventory be stored in containers. From a
mathematical perspective, containers are  elements of a finite or countably infinite partition of $\X=\mathbb{R}$ into disjoint convex sets, and each of
these sets is not a singleton. In other words, each
 container $B_{i+1}$  is an interval (possibly open, closed, or
semi-open) with ends $d_i$ and $d_{i+1}$ such that $-\infty\le
d_i<d_{i+1}\le +\infty$, and the union of these disjoint intervals
is $\mathbb{R}.$ In addition, we assume that
$d_{i+1}-d_i\ge\gamma$ for some constant $\gamma>0$ for all
containers, that is, the sizes of all the containers are uniformly
bounded below by a positive number.  We also follow the convention
that the 0-inventory level belongs to a container with end points
$d_0$ and $d_1$, and a container with end points $d_i$ and
$d_{i+1}$ is labeled as the $(i+1)$-th container $B_{i+1}$. Thus,
container $B_1$ is the interval in the partition containing point
0. The containers' labels can be nonpositive.  If there is a container
with the smallest (or largest) finite label $n$ then
$d_{n-1}=-\infty$ (or $d_n=+\infty$, respectively).  If there are
containers with  labels $i$ and $j$ then there are containers with
all the labels between $i$ and $j$. In addition each container is
either transparent or nontransparent. If the inventory level $x_t$
belongs to a nontransparent container, the IM only  knows which
container the inventory level belongs to. If an inventory level
$x_t$ belongs to a transparent container, the IM knows that the
amount of inventory is exactly $x_t;$ see Figures~\ref{inventory_knownLevel}--\ref{inventory_unknownZeroLv}.

For each nontransparent container with end points $d_i$ and
$d_{i+1}$, we fix an arbitrary point $b_{i+1}$ satisfying
$d_i<b_{i+1}<d_{i+1}$. For example, it is possible to set
$b_{i+1}=0.5d_i+0.5d_{i+1},$ when $\max\{|d_i|, |d_{i+1}|\}<\infty.$
If an inventory level belongs to a nontransparent container $B_{i}$,
the IM observes $y_t=b_i.$ Let $L$ be the set of labels of the
nontransparent containers.  We set $Y_L=\{b_i\,:\,i\in L\}$ and
define the observation set $\Y=\T\cup Y_L$, where $\T$ is the union
of all transparent containers $B_i$ (transparent elements of the
partition).  If the observation $y_t$ belongs to a transparent
container (in this case, $y_t\in\T$), then the IM knows that the
inventory level $x_t=y_t$. If $y_t\in Y_L$ (in this case, $y_t=b_i$
for some $i$), then the IM knows that the inventory
level belongs to the container $B_i$, 
and this container is nontransparent. Of course, the distribution
of this level can be computed. 

Let $\rho$ be the Euclidean distance on $\mathbb{R}:$ $\rho
(a,b)=|a-b|$ for $a,b\in \Y$. On the state space $\X=\mathbb{R}$ we
consider the metric $\rho_\X(a,b)=|a-b|,$ if $a$ and $b$ belong to
the same container, and $\rho_\X(a,b)=|a-b|+1$ otherwise, where
$a,b\in \X$. The space $(\X,\rho_\X)$ is a Borel subset of a Polish
space (consisting of closed containers, that is, each finite point
$d_i$ is represented by two points: one belonging to the container
$B_i$ and another one to the container $B_{i+1}$).
 We notice that
 $\rho_\X(x^{(n)},x)\to 0$ as $n\to\infty$  if and only if $|x^{(n)}-x|\to
0$ as $n\to\infty$ and the sequence $\{x^{(n)}\}_{n=N,N+1,\ldots}$
belongs to the same container as $x$ for a sufficiently large $N$.
Thus, convergence on $\X$ in the metric $\rho_{\X}$ implies
convergence in the Euclidean metric.  In addition, if $x\ne d_i$
for all containers $i$, then $\rho_\X(x^{(n)},x)\to 0$ as
$n\to\infty$ if and only if $|x^{(n)}-x|\to 0$ as $n\to\infty.$
Therefore, for any open set $B$ in $(\X,\rho_\X)$, the set
$B\setminus(\cup_i\{d_i\})$ is open in $(\X,\rho).$  We notice
that each container $B_i$ is an open and closed set in
$(\X,\rho_\X).$

It is possible to show that the state transition law $P$ given by \eqref{eq:P} is weakly continuous in $(x,a)\in\X\times\A$.
Set $\Psi(x)=x,$ if the inventory level $x$ belongs to a transparent container, and $\Psi(x)=b_i,$ if the inventory level belongs to a nontransparent
container $B_i$ with a label $i$. As follows from the definition of the metric $\rho_\X$, the function $\Psi:(\X,\rho_\X)\to (\Y,\rho)$ is continuous.
Therefore, the observation transition probabilities  $Q_0$   from $\X$  to $\Y$ and $Q$ from $\A\times\X$   to $\Y ,$ $Q_0(C|x):=Q(C|a,x):=\h\{\Psi(x)\in C\}$, $C\in \B(\Y)$,
$a\in\A$, $x\in\X$, are weakly continuous.

If all the containers are nontransparent, the observation set
$\Y=Y_L$ is countable,  and conditions of Corollary~\ref{cor:main1}
hold. In particular, the function $Q(b_i|a,x)=\h\{x\in B_i\}$ is
continuous, if the metric $\rho_\X$ is considered on $\X.$ If some
containers are transparent and some are not, the conditions of
Corollary~\ref{cor:main1} hold. To verify this, we set $\Y_1:=\T$
and $\Y_2:=Y_L$ and note that $\Y_2$ is countable  and the function
$Q(b_i|x)=1\{x\in B_i\}$ is continuous for each $b_i\in Y_L$
because $B_i$ is open and closed in $(\X,\rho_\X).$ Note that
$H(B|z,a,y)=P( B|y,a)$ for any $B\in \B(\X)$, $C\in\B(\Y)$,
$z\in\P(\X)$, $a\in \A$, and $y\in\T$. The kernel $H$ is weakly
continuous on $\P(\X)\times\A\times \Y_1$. In addition, $\T=\cup_i
B^{tr}_i$, where $B_i^{tr}$ are transparent containers, is an open set in
$(\X,\rho_\X).$ Thus the POMDP ($\X$, $\Y$, $\A$, $P$, $Q$, $c$)
satisfies the assumptions of Corollary~\ref{cor:main1}. Thus, for
the corresponding COMDP, there are stationary optimal policies,
optimal policies satisfy the optimality equations, and value
iterations converge to the optimal value.

The models studied in Bensoussan et al. \cite{Bens3, Bens4, Bens5} correspond to the partition $B_1=(-\infty,0]$ and $B_2=(0,+\infty)$   with the container
$B_2$ being nontransparent and with the container $B_1$ being either nontransparent (backordered amounts are not known \cite{Bens5}) or transparent (models
with lost sales \cite{Bens4},
 backorders are observable \cite{Bens3}). 
Note that, since $F_D$ is atomless, the probability that
$x_t+a_t-D_{t+1}=0$ is $0$, $t=0,1,\ldots\ .$

The model provided in this subsection is applicable to other inventory control problems, and the conclusions of Corollary~\ref{cor:main1} hold for them
too.   For example, consider a periodic review inventory system with backorders, for which nonnegative inventory levels are known, and, when the inventory level is negative,  it is known that
there is a backorder, but its quantity is unknown.  The partition consists of two containers: a nontransparent container $B_{0}=(-\infty,0)$ and a transparent container $B_1=[0, +\infty).$
%
%

\section{Conclusions}
The tutorial describes general sufficient conditions for the existence and characterization of optimal policies for Markov Decision Processes with possibly infinite state spaces and unbounded action sets and costs.  Expected total discounted cost and average cost criteria are considered.  The described conditions imply the existence of optimal Markov policies in finite-horizon problems and the existence of optimal stationary policies for infinite-horizon problems.  They imply the validity of optimality equations, convergence of value iterations, and continuity properties of value functions for discounted costs.  They also imply the validity of optimality inequalities for average costs per unit time.

For discounted costs, these conditions consist of two assumptions: the transition probabilities are weakly continuous, and the one-step cost function is $\K$-inf-compact.  These two assumptions practically always hold for periodic-review stochastic inventory control problems.  The $\K$-inf-compactness property of one-step costs is weaker than inf-compactness, which typically holds for cost functions for inventory control problems.  One of the reasons for the generality of the results is that their derivation is linked to a new maximum theorem, which extends Berge's maximum theorem to possibly noncompact action sets.

For average cost MDPs, the single additional assumption is  that the relative value function is well-defined. This assumption also holds for inventory control applications and can be verified  easily.

The tutorial also describes optimality conditions for Partially Observable Markov Decision Processes with total discounted costs.  These conditions imply the existence of optimal policies, validity of optimality equations, and convergence of value iterations.  The results are illustrated with inventory control models for which some of the inventory levels are not observable.

 The described results and methods are useful and insightful for
 investigating new and existing inventory control problems.  As an
 illustration, a complete classification of possible solutions for the
 classic periodic-review stochastic single-product problem is
 described.

\section*{Appendix}

\section{Berge's Maximum Theorem for Noncompact Action Sets and Some Properties of $\K$-Inf-Compact Functions}\label{Apx}
This appendix describes generalizations of Berge's maximum theorem and the relevant Berge
theorem on semi-continuity of the value function to possibly noncompact action sets.
These theorems are important for control theory, games, and mathematical economics.  The major limitation of these  theorems is that they require compact action sets. The generalizations provided in Feinberg et al., \cite{FKV, FKZ} remove this limitation.  Here we present these results for metric spaces.  With slight modifications they hold for Hausdorff topological spaces
(see  \cite{FKV}), but this level of generality is not needed for the results of this tutorial. Local versions of the results presented in this appendix can be found in Feinberg and Kasyanov~\cite{FK}.

Let ${\Sp^1}$ and ${\Sp^2}$ be
metric spaces, $u:{\Sp^1}\times {\Sp^2}\to
\overline{\mathbb{R}}=\mathbb{R}\cup\{\pm\infty\}$ and
$\Phi:{\Sp^1}\to 2^{{\Sp^2}}\setminus \{\emptyset\}$. Consider an
optimization problem of the form
\begin{equation}\label{eq:lm10}
v({s^1}):=\inf\limits_{{s^2}\in {\Phi}({s^1})}u({s^1},{s^2})\quad\mbox{for each}\quad
{s^1}\in{\Sp^1};
\end{equation}
which appears, for instance, in optimal control and game theory.
  Let 
  $\K({\Sp^2})$ be the set of nonempty compact subsets of ${\Sp^2}.$  Berge's theorem has the following formulation.\\

{\bf Berge's Theorem.} (\cite[p. 116]{Ber}). If
$u:{\Sp^1}\times{\Sp^2}\to\overline{\mathbb{R}}$ is a lower semi-continuous
function and $\Phi:{\Sp^1}\to \K({\Sp^2})$ is an
upper semi-continuous set-valued mapping, then the
function $v:{\Sp^1}\to\overline{\mathbb{R}}$ is lower semi-continuous.
\vskip 0.2 cm

The well-known Berge's  maximum
theorem
 has the following formulation.\\
{\bf Berge's Maximum Theorem.} (\cite[p. 116]{Ber}). If
$u:{\Sp^1}\times{\Sp^2}\to{\mathbb{R}}$ is a continuous function and
$\Phi:{\Sp^1}\to \K({\Sp^2})$ is a 
continuous set-valued
mapping, then the value function $v:{\Sp^1}\to\R$ is continuous and
the solution multifunction ${\Phi}^*:{\Sp^1}\to 2^{{\Sp^2}}\setminus \{\emptyset\}$, defined as
\begin{equation}\label{eq:lm000} {\Phi}^*({s^1})=\left\{{s^2}\in
{\Phi}({s^1}):\,v({s^1})=u({s^1},{s^2})\right\},\quad {s^1}\in{\Sp^1},
\end{equation}
is  upper semi-continuous and compact-valued.
\vskip 0.2 cm

 For an $\overline{\mathbb{R}}$-valued
function $f$, defined on a nonempty subset $U$ of a topological
space $\mathbb{U},$ consider the level sets
\begin{equation*}
\mathcal{D}_f(\lambda;U)=\{y\in U \, : \,  f(y)\le
\lambda\},\qquad \lambda\in\R.\end{equation*}  We recall that a
function $f$ is \textit{lower semi-continuous on $U$} if all the
level sets $\mathcal{D}_f(\lambda;U)$
 are closed, and a function $f$ is
\textit{inf-compact} (also sometimes called \emph{lower
semi-compact}) on $U$ if all these sets are compact. The following definition deals
with the space $\U={\Sp^1}\times{\Sp^2}$ and its subsets ${\rm Gr}_{{\Sp^1}}(\Phi)$ and ${\rm Gr}_K(\Phi).$

\begin{defn}\label{DKKFUN} (\cite[Definition 1.1]{FKZ}).
A function $u:{\Sp^1}\times {\Sp^2}\to \overline{\mathbb{R}}$ is called
$\K$-inf-compact on ${\rm Gr}_{{\Sp^1}}(\Phi)$, 
if for every compact subset $K$ of ${\Sp^1}$ this function is inf-compact on ${\rm
Gr}_K(\Phi)$.
\end{defn}
 The following two theorems generalize  Berge's theorem and Berge's maximum theorem respectively to possibly noncompact action sets.
\begin{thm}\label{MT1} (\cite[Theorem 1.2]{FKZ}).
If the function $u:{\Sp^1}\times {\Sp^2}\to \overline{\mathbb{R}}$ is
$\K$-inf-compact on ${\rm Gr}_{{\Sp^1}}(\Phi)$, then the function
$v:{\Sp^1}\to\overline{\mathbb{R}}$ is lower semi-continuous.
\end{thm}
\begin{thm}\label{thm: Berge}(\cite[Theorem 1.2]{FKV}).
Assume that:
\begin{enumerate}[(a)]
%
\item $\Phi: {\Sp^1} \to 2^{\Sp^2}\setminus\{\emptyset\}$ is lower semi-continuous; \item $u: {\Sp^1}
\times {\Sp^2} \to \R$ is $\K$-inf-compact and upper semi-continuous on
${\rm
Gr}_{{\Sp^1}}(\Phi)$.
\end{enumerate}
Then the value function $v: {\Sp^1}\to \R$ is continuous and the solution multifunction $\Phi^*:{\Sp^1}\to \K({\Sp^2})$ is 
upper semi-continuous and compact-valued.
\end{thm}

The first statement of the following lemma implies that Theorems~\ref{MT1} and \ref{thm: Berge} are indeed generalizations of Berge's theorem and Berge's maximum theorem respectively.  The second statement indicates that the class of $\K$-inf-compact functions is broader
than the class of inf-compact functions.

\begin{lem}\label{lm0} (\cite[Lemma 2.1]{FKZ}).
The following statements hold:

(i) if $u:{\Sp^1}\times{\Sp^2}\to \overline{\mathbb{R}}$ is lower
semi-continuous on ${\rm Gr}_{\Sp^1} ({{\Phi}})$ and ${\Phi}:{\Sp^1}\to
\K({\Sp^2})$ is upper semi-continuous, then the function
$u(\cdot,\cdot)$ is $\K$-inf-compact on ${\rm Gr}_{\Sp^1} ({{\Phi}})$;

(ii) if $u:{\Sp^1}\times{\Sp^2}\to \overline{\mathbb{R}}$ is inf-compact on
${\rm Gr}_{\Sp^1}(\Phi)$, then the function $u(\cdot,\cdot)$ is
$\K$-inf-compact on ${\rm Gr}_{\Sp^1}(\Phi)$.
\end{lem}

Luque-V\'asquez and Hern\'andez-Lerma \cite{LVHL} provided an
example with ${\Sp^1}=\R,$ ${\Sp^2}=\Phi({s^1})=[0,\infty),$  continuous $\Phi,$ and continuous $u({s^1},{s^2})$  which is inf-compact in ${s^2},$ where $v({s^1})$ is not lower semi-continuous.  The following two lemmas indicate that $\K$-inf-compactness of $u$ is stronger than its lower-semicontinuity and inf-compactness in ${s^2}.$

\begin{lem}\label{lm00101} (\cite[Lemma 2.2]{FKZ}).
If $u(\cdot,\cdot)$ is $\K$-inf-compact function on ${\rm
Gr}_{{\Sp^1}}(\Phi)$, then for every ${s^1}\in{\Sp^1}$ the function $u({s^1},\cdot)$
is inf-compact on ${\Phi}({s^1})$.
\end{lem}

\begin{lem}\label{lsc} (\cite[Lemma 2.3]{FKZ}).
A $\K$-inf-compact function $u(\cdot,\cdot)$ on ${\rm
Gr}_{{\Sp^1}}(\Phi)$ is lower semi-continuous on ${\rm Gr}_{{\Sp^1}}(\Phi)$.
\end{lem}

The following lemma provides the necessary and sufficient condition for $\K$-inf-compactness.
This condition is used in Assumption~{\bf W*} in Feinberg et al,~\cite{FKZMDP} instead of equivalent Definition A.1.
\begin{lem}\label{lm4} (\cite[Lemma 2.5]{FKZ}).
The function $u(\cdot,\cdot)$ is
$\K$-inf-compact on ${\rm Gr}_{{\Sp^1}}(\Phi)$ if and only if the
following two conditions hold:

(i) $u(\cdot,\cdot)$ is lower semi-continuous on ${\rm
Gr}_{{\Sp^1}}(\Phi)$;

(ii) if a sequence $\{s^1_n \}_{n=1,2,\ldots}$ with values in ${\Sp^1}$
converges and its limit $s^1$ belongs to ${\Sp^1}$ then any sequence
$\{s^2_n \}_{n=1,2,\ldots}$ with $s^2_n\in \Phi(s^1_n)$, $n=1,2,\ldots,$
satisfying the condition that the sequence $\{u(s^1_n,s^2_n)
\}_{n=1,2,\ldots}$ is bounded above, has a limit point $s^2\in
\Phi(s^1).$
\end{lem}
\textbf{Acknowledgement} Some of the materials presented in this tutorial are based on results of work partially supported by NSF grant CMMI-1335296.  The author thanks Jefferson Huang, Pavlo O. Kasianov, Mark E. Lewis, Yan Liang, and Matthew J. Sobel for valuable comments.





\end{document}